\documentclass[lettersize,journal]{IEEEtran}
\usepackage{amsmath,amsfonts}
\usepackage{algorithmic}
\usepackage{algorithm}
\usepackage{array}
\usepackage[caption=false,font=normalsize,labelfont=sf,textfont=sf]{subfig}
\usepackage{textcomp}
\usepackage{stfloats}
\usepackage{url}
\usepackage{verbatim}

\usepackage{graphicx} 
\usepackage{color,soul}
\usepackage[utf8]{inputenc}

\usepackage{amsthm}
\usepackage{amssymb}
\usepackage{esvect}

\DeclareMathAlphabet{\mathbit}{OML}{cmm}{b}{it}
\usepackage{mathtools}
\usepackage{accents}

\usepackage{stackengine}
\usepackage{comment}
\usepackage[version=4]{mhchem}
\usepackage{siunitx}
\usepackage{longtable,tabularx}
\DeclarePairedDelimiter\norm{\lVert}{\rVert}%
\usepackage{pdfpages}
\usepackage{mathtools,lipsum}

\usepackage{cuted}
\usepackage{bm}

\begin{document}
\bstctlcite{bstctl:nodash}

\title{On The Dynamics Of Variable-Shape Wave Energy Converters}

\author{Mohamed A. Shabara, and Ossama Abdelkhalik
\thanks{Mohamed A. Shabara and Ossama Abdelkhalik are with the Department of Aerospace Engineering, Iowa State University, Ames, IA, USA, 50010 (e-mail: mshabara@iastate.edu;
ossama@iastate.edu.com)}}



\maketitle

\begin{abstract}
Flexible structures in the wave energy conversion field have recently captured attention. The common method used to tackle the problem of flexibility in wave energy devices is done by accounting for the flexibility (mode shapes) in the hydrodynamic solvers as extra degrees of freedom, not in the structural domain. This work uses the principles of analytical mechanics to derive the equation of motion of spherical variable-shape wave energy converters and describes a methodology to calculate the generalized hydrodynamic forces on the shell. Fluid-Structure Interaction simulations are performed to validate the developed dynamic model and to study the effect of using a flexible buoy in wave energy converters on its trajectory and power production. 
\end{abstract}


\section{Introduction}
\IEEEPARstart{O}{ne} of the early studies on the vibration of thin elastic shells was carried out by Love \cite{love1887small}, which resulted in what is known as Love's approximation. The formulation in \cite{love1887small} laid the foundation for the classical bending theory \cite{naghdi1962vibrations}. Spherical shells vibration has multiple applications such as spherical pressure vessels/tanks, space vehicles like the tumbleweed rovers used in Mars explorations \cite{9144386,forbes2010dynamic,wilson2008design}, dynamics of gas bubbles, bio-medical applications, naval vehicles, and wave energy converters ``WECs" \cite{9207738}. 


The wave energy converters can be classified either based on their position (shoreline, near-shore and offshore devices), or principle of operation, or based on the power take-off (PTO) technique. Drew et al. \cite{drew2009review} classified the WEC based on the principle of operation into three main categories as follows: \begin{enumerate}
\item Attenuators, such as the Pelamis are devices that lie parallel to the predominant wave direction.
\item Point absorbers, have small dimensions relative to the incident wave wavelength.  
\item Terminators, such as the Salter's duck are devices that lie perpendicular to the predominant wave direction.
\end{enumerate} 

Typical optimal control strategies for WECs require reactive power; reactive power flows from the device to water and usually needs a relatively complex Power Take-off (PTO) units. Another aspect is the power quality (fluctuations in the power curve) \cite{zou2020modeling,zou2021numerical}. On the other hand, a passive control does not need reactive power; yet the converted power is not as high as that obtained when using reactive power.

Variable-shape-buoy wave energy converters (VSB WECs) were recently introduced to tackle this reactive power challenge by either reducing the need for it or completely eliminating the reactive power component. It is worth noting that the difference between a VSB WEC and a variable geometry WEC is the rate of change of the surface of the WEC shape relative to the incident wave forces, i.e., the variable geometry WEC is a Fixed Shape buoy (FSB) WEC that changes its geometry occasionally, on the other hand, the VSB WEC changes their shapes continuously. The wave/WEC interaction produced by the VSB WEC can be tuned to produce more power without the complexity accompanied by the reactive power PTO units. 

Flexible Wave Energy Converters (FlexWecs) uses the Distributed embedded energy conversion technology (DEEC-Tec) to harvest energy from the ocean, the structure of these devices use distributed smart materials across their shells to harvest wave energy from the waves \cite{boren2021distributed,boren2022wbs}. It is worth noting that the work in this paper can be extended to model such devices by including the piezoelectric effects in the model.

Multiple optimal control strategies were introduced in the literature; among the most implemented passive control methods are the optimal resistive control \cite{zou2017optimal}, passive loading \cite{garcia2017real}, latching control \cite{budar1975resonant}, and passive MPC \cite{montoya2021passive}. 
The reactive power capability enables the PTO to not only harvest energy but also derive the buoy at certain times to create resonance in the WEC system (mechanical impedance matching). Optimal control methods derived using Pontryagin’s minimum principle result in reactive powers \cite{zou2017optimal}. Habeebullah et al. \cite{9705961} applied the optimal control theory to obtain an optimal control solution for the PTO while constraining the reactive power not to exceed a certain threshold, i.e., reducing the complexity of the PTOs while increasing the harvested energy compared to passive control methods.    

Zou et al. \cite{zou2020modeling} proposed a design for a VSB WEC that is composed of a cylindrical gas chamber with $2$ m radius attached from the bottom to a set of multiple controllable/movable panels. The latter VSB WEC was controlled uisng a simple linear damping power take-off unit ``PTO". The device interacts non-linearly with the incident waves (including deformation and translation). A low-fidelity dynamic model is derived to validate the superiority of the VSB WEC over the FSB WEC. The power harvesting using the VSB WEC was almost 18\% higher when compared to the FSB WEC. Also, it is noticed that the $pk$-$pk$ velocity of the VSB WEC is higher compared to the FSB WEC.

A Fluid-Structure Interaction (FSI) high fidelity simulation was conducted (using computational fluid dynamics and finite element analysis) in \cite{shabara2021numerical} to assess the performance of a spherical VSB WEC compared to a FSB WEC. A concrete plate was attached to the middle section of the VSB WEC that divides the internal volume into two separate partitions. A passive control approach is adopted; the control force is applied to the concrete plate. The results showed an increase in the heave displacement and velocity responses; the increase in $pk-pk$ heave response was $11.88\%$. The results also showed an increase in the harvested energy of 90 KJ over a time interval of $30$ seconds. It is worth noting that this increase was mainly due to the transient effects in the first few seconds of the simulation.

The most commonly used equation of motion for WECs is Cummin's Equation which is expressed as: \begin{align}\label{EoMF}
m\ddot{z}(t) = &\overbrace{\int_{-\infty}^{\infty}h_f(\tau)\eta(t - \tau,z)d\tau}^\text{excitation force $f_e$} + f_s \nonumber\\  & \qquad \qquad  \overbrace{-\mu\ddot{z}(t) - \int_{-\infty}^{t}h_r(\tau)\dot{z}(t - \tau)d\tau}^\text{radiation force $f_r$} - u
\end{align} where $m$ is the mass of the buoy, $z$, $\dot{z}$ and $\ddot{z}$ are the heave, velocity and acceleration of the buoy's center of gravity (C.G). $f_e$ is the excitation force due to the impulses $h_f(\tau) $, $\eta$ is the significant wave height. $f_r$ is the force applied to the buoy due to the radiated waves, which it is dependent on both the velocity and the acceleration of the device and $\mu$ is the added mass. The convolution term in the $f_r$ in the state space model as: \begin{eqnarray}\label{eqXr}
\dot{\vec{x}}_r  =  A_r \vec{x}_r + B_r \dot{z} & \text{and} & f_r  =  C_r \vec{x}_r,
\end{eqnarray} where $A_r$, $B_r$ and $C_r$ are derived based on the impulse response function in the Laplace domain. 

Eq.~\eqref{EoMF} can be further extended to include the six modes of motion (degrees of freedom) for FSB WECs \cite{bacelli2014optimal}; however, it is not used to solve for VSB WECs as it does not account for flexibility. Different approaches were developed to account for the flexibility of VSB WECs, however, these methods either account for the flexibility modes of the WEC in the fluid solver \cite{mcdonald2019linear}, or was limited only to simple geometry WECs \cite{guo2017inclusion}.


This work develops a rigorous dynamic model for VSB WECs based on the fundamentals of analytical mechanics and shell theories. The paper derives the motion equation for spherical VSB WECs. The presented approach can be applied for other shapes of VSB WECs. The main concern in this work is to understand the structural and dynamic behavior of the VSB WEC shell; accordingly, a relatively simple sinusoidal wave is assumed to act on the shell. This paper is divided into 6 main sections. In section \ref{sec:math}, the kinetic and potential energies of asymmetric free vibrating spherical shells is derived using a similar approach to the approach presented in \cite{forbes2010dynamic,hogan2015dynamic,robert2016modeling}. In section \ref{sec:EOM}, the equations of motion for spherical shell buoys are derived using Lagrangian mechanics for the free unconstrained case. The equation of motion for the forced case is derived in section \ref{sec:const_forced}. The numerical simulation results are discussed in section \ref{sec:results}.   

\section{Kinetic and Potential Energies of Spherical Shell Buoys} \label{sec:math}

Hamilton principle of least action is based on the minimization of system's energy, which is the base of the Lagrangian mechanics approach implemented in this work; therefore, the knowledge of the total energy of the system is required (i.e., the total kinetic and total potential energies should be computed). To compute the total energy of the VSB system a kinematic analysis is carried out which starts with defining the reference frames (subsection \ref{subsec:Kinematics}), them using these reference frames to write expressions for the kinetic energy (subsection \ref{subsec:KE}) and potential energy (subsections \ref{subsec:PE}).

\begin{figure}
\centering
\includegraphics[scale=0.65]{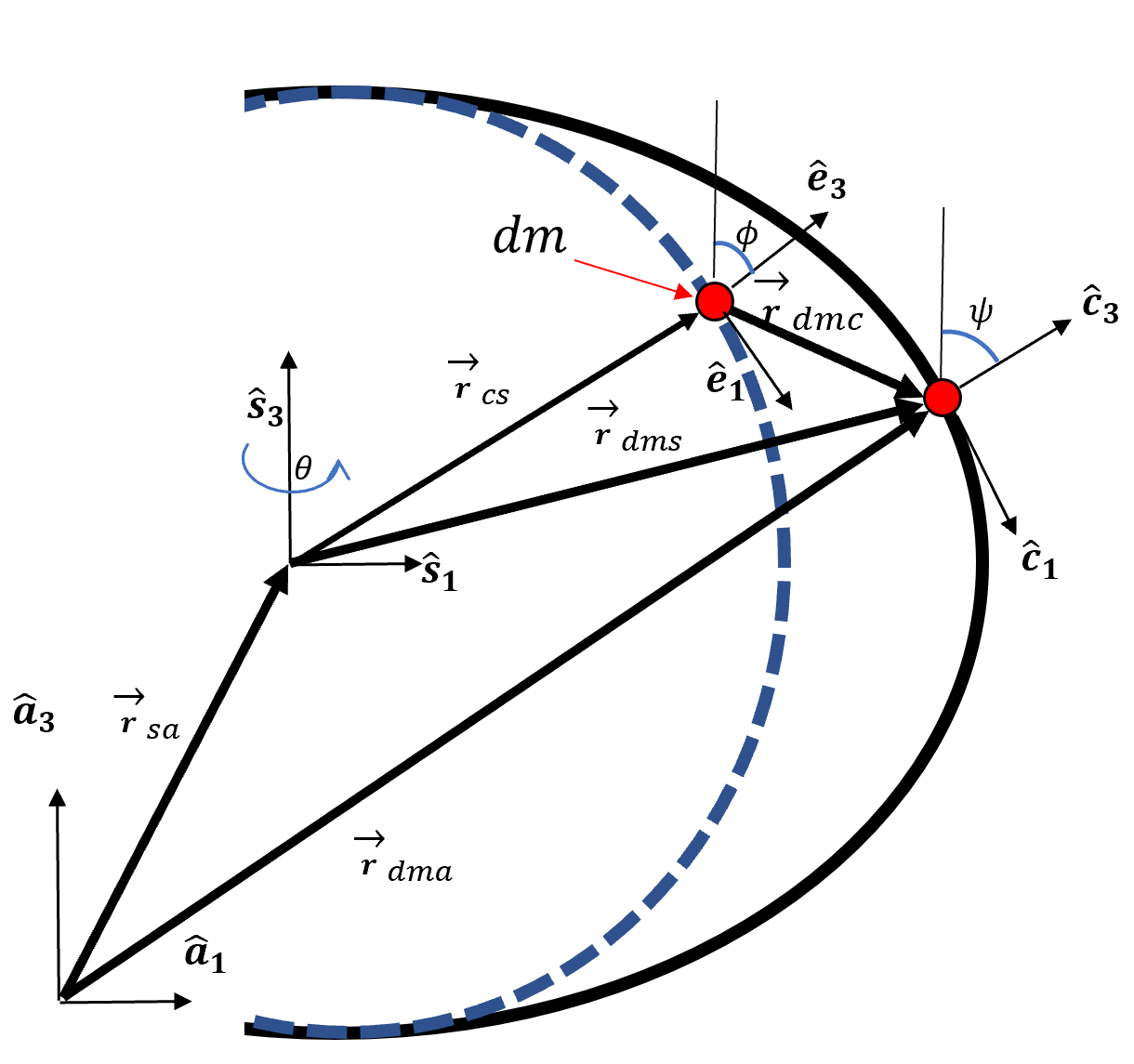}
\caption{Deformed (Solid Black Line) and non-deformed Sphere (Dashed Blue Line)}
\label{fig:1}
\end{figure}

Consider a VSB for which the non-deformed shape (not necessarily spherical) is shown Fig. (\ref{fig:1}), the inertial frame ($\hat{\boldsymbol{a}}$) can be expressed as:
\begin{equation*}
\hat{\boldsymbol{a}} = \left[ \hat{\boldsymbol{a}}_1 , \hat{\boldsymbol{a}}_2 , \hat{\boldsymbol{a}}_3 \right ] 
\end{equation*}

Consider a body-fixed frame $\hat{\boldsymbol{s}}$ that is attached to the buoy's C.G. 
The location of any infinitesimal mass on the buoy's shell can be specified using the two angels $\phi$ and $\theta$, as shown in Fig. (\ref{fig:1}).
Consider a the reference frame $\hat{\boldsymbol{e}}$ attached to an infinitesimal mass ($dm$) at the surface of the buoy before deformation, and its third axis $\hat{\boldsymbol{e}}_3$ points at the radial direction of the non-deformed buoy shape.
Hence, the reference frame $\hat{\boldsymbol{e}}$ is obtained by rotating $\hat{\boldsymbol{s}}$ through an angle $\theta$ around the $\hat{\boldsymbol{s}}_3$ then through an angle $\phi$ about the second axis of the intermediate frame as follows: 
\begin{equation*}
    C_{es}(\phi,\theta) = C_2(\phi)C_3(\theta) 
\end{equation*}where $C_i(x)$ represents a cosine transformation matrix of a single rotation of angle $x$ about the coordinate $i$, $i=1,2$ and $3$.
The reference frame $\hat{\boldsymbol{c}}$ is attached to the infinitesimal mass ($dm$) on the surface after the shell deformation such that $\hat{\boldsymbol{c}}_3$ is normal to the shell surface.

The angle between $\hat{\boldsymbol{s}}_3$ and $\hat{\boldsymbol{c}}_3$ axes is $\psi$. In this analysis, it is assumed that the deformations are axisymmetric about the $\hat{\boldsymbol{a}}_3$ axis; thus, the axis $\hat{\boldsymbol{c}}_2$ is always perpendicular to the page so as the axes $\hat{\boldsymbol{a}}_2$, $\hat{\boldsymbol{s}}_2$, and $\hat{\boldsymbol{e}}_2$.
For a non deformed shell, the frames $\hat{\boldsymbol{e}}$ and $\hat{\boldsymbol{c}}$ coincide. The frames $\hat{\boldsymbol{e}}$ and $\hat{\boldsymbol{c}}$ become different, in general, when the shape is deformed. 
For the FSB WEC the reference frames $\hat{\boldsymbol{e}}$ and $\hat{\boldsymbol{c}}$ coincide; this applies for the VSB WEC at the initial time before deformation.

Noting that, high fidelity simulations were carried out in References \cite{zou2021numerical,shabara2021numerical,zou2020numerical} and it was found that the VSB deform in a nearly axisymmetric behaviour at the steady state response. i.e., the axisymmetric vibration assumption in this work.  

The coordinate transformation matrix from the $\hat{\boldsymbol{a}}$ frame to the $\hat{\boldsymbol{s}}$ frame is computed in this paper using the the 3-2-1 Euler angle sequence  as $C_{sa}(\alpha,\beta,\gamma) =  C_1(\alpha)C_2(\beta) C_3(\gamma)$. The cross product of two arbitrary vectors expressed in the same reference frame can be replaced with matrix multiplication:  $\vec{\boldsymbol{u}} \times \vec{\boldsymbol{v}} = \boldsymbol{u}^\times \vec{\boldsymbol{v}}$, where \begin{equation*}
\boldsymbol{u}_a^\times = \begin{bmatrix}
0 & -u_{a3} & u_{a2}\\
u_{a3} & 0 & -u_{a1} \\
-u_{a2} & u_{a1} & 0  \\ 
\end{bmatrix} \\
\end{equation*} Since the changes in the buoy shape are assumed axisymmetric, we can express the deformation vector (displacement) ``$ \vec{\boldsymbol{r}}_{dmc}(\theta,\phi,t)$" as function of only the angle $\phi$ and the time $t$. This deformation vector can be expressed in the $\hat{\boldsymbol{e}}$ frame as \cite{robert2016modeling,hogan2015dynamic}: \begin{equation} \label{eq:1}
 \vec{\boldsymbol{r}}_{dmc}(\phi,t) = \begin{bmatrix} u (\phi,t) & 0 & v (\phi,t) \end{bmatrix}{}^T
\end{equation} where the second component ${{r}_{dmc}}_2$ is set to zero due to the axisymmetry assumption, $u(\phi,t)$ is the displacement component in the $\hat{\boldsymbol{e}}_1$ direction, and $v(\phi,t)$ is the displacement component in the $\hat{\boldsymbol{e}}_3$ direction. Solving the problem using the above expression for the deformation vector yields a distributed parameter model which is rather complicated to solve as illustrated in the following section, thus, each term of the deformation vector is assumed to be a series of separable functions and Rayleigh-Ritz approximation is summoned. 


\subsection{Kinematics of a Flexible Spherical Buoy -  Free Vibration} \label{subsec:Kinematics}
This subsection is concerned with calculating the velocity vector ${}^a\dot{\vec{\boldsymbol{r}}}_{dma}$ as it is used in the calculation of the kinetic energy of the shell due to the translation and rotational motions as well as the deformation of the external shell.
As shown in Fig. (\ref{fig:1}), the position vector of a point on the surface of the deformed sphere in the inertial frame ``$\hat{\boldsymbol{a}}$" is expressed as \begin{equation} \label{eq:24}
\vec{\boldsymbol{r}}_{dma}=\vec{\boldsymbol{r}}_{sa} +\vec{\boldsymbol{r}}_{cs} +\vec{\boldsymbol{r}}_{dmc}
\end{equation}
The velocity vector is expressed as \cite{robert2016modeling, hogan2015dynamic}:
\begin{equation}
{}^a\dot{\vec{\boldsymbol{r}}}_{dma} ={}^a\dot{\vec{\boldsymbol{r}}}_{sa} +{}^a\dot{\vec{\boldsymbol{r}}}_{cs} +{}^a\dot{\vec{\boldsymbol{r}}}_{dmc} 
\end{equation} Note that the left superscript denotes the reference frame used to describe the vector. Applying the transport theorem knowing that ${}^s\dot{\vec{\boldsymbol{r}}}_{cs} = \boldsymbol{0}$, we get:
\begin{align}
{}^a\dot{\vec{\boldsymbol{r}}}_{dma} & = {}^a\dot{\vec{\boldsymbol{r}}}_{sa} +   ({}^s\dot{\vec{\boldsymbol{r}}}_{cs} + \vec{\omega}_{as}\times{}^s\vec{\boldsymbol{r}}_{cs}) \nonumber \\ & \qquad \quad \; + ({}^s\dot{\vec{\boldsymbol{r}}}_{dmc} +\vec{\omega}_{as} \times {}^s\vec{\boldsymbol{r}}_{dmc}  \ )\\
&={}^a\dot{\vec{\boldsymbol{r}}}_{sa} +{}^s\dot{\vec{\boldsymbol{r}}}_{dmc} -( {}^s\vec{\boldsymbol{r}}_{cs} + {}^s\vec{\boldsymbol{r}}_{dmc}) \times  \vec{\omega}_{as} \label{eq:27}
\end{align}

Note that ${}^e\vec{\boldsymbol{r}}_{cs} = r \hat{\boldsymbol{e}}_3$; hence by substituting in Eq. (\ref{eq:27}) we obtain: \begin{align}
{}^a\dot{\vec{\boldsymbol{r}}}_{dma} = {}^a\dot{\vec{\boldsymbol{r}}}_{sa} &+ {}^e\dot{\vec{\boldsymbol{r}}}_{dmc} + \vec{\omega}_{es} \times \vec{r}_{dmc}
\nonumber \\ &- (C_{se} {}^e\vec{\boldsymbol{r}}_{cs} + C_{se}  {}^e\vec{\boldsymbol{r}}_{dmc}) \times  \vec{\omega}_{as} 
\end{align}

Note that $\vec{\omega}_{es} = 0$, and express the first two terms in the above equation in the $\hat{s}$ frame to get:
\begin{equation} 
{}^a\dot{\vec{\boldsymbol{r}}}_{dma}  = C_{sa}{}^a\dot{\vec{\boldsymbol{r}}}_{sa} + C_{se} {}^e\dot{\vec{\boldsymbol{r}}}_{dmc}
- [C_{se} ({}^e\vec{\boldsymbol{r}}_{cs} 
+ {}^e\vec{\boldsymbol{r}}_{dmc})]^\times  \vec{\omega}_{as} \end{equation}
\begin{align} & = \begin{bmatrix}
C_{sa}  & -[C_{se} ({}^e\vec{\boldsymbol{r}}_{cs} + {}^e\vec{\boldsymbol{r}}_{dmc})]^\times  & C_{se} 
\end{bmatrix} \underbrace{\begin{bmatrix}
\vec{\omega}_{as} \\
{}^e\dot{\vec{\boldsymbol{r}}}_{dmc}
\end{bmatrix}}_{\dot{\vec{\mathbf{x}}}} \\
&= \begin{bmatrix}
C_{sa}  & -[C_{se} ({}^e\vec{\boldsymbol{r}}_{cs} + {}^e\vec{\boldsymbol{r}}_{dmc})]^\times  & C_{se} 
\end{bmatrix} \dot{\vec{\mathbf{x}}} \label{eq:31}
\end{align}

\subsection{Kinetic Energy for a Flexible Spherical Buoy } \label{subsec:KE} 
The total kinetic energy of the buoy is expressed in Eq. (\ref{eq:32}).
\begin{equation} \label{eq:32}
\mathcal{T} = \frac{1}{2} \int_S   {}^a\dot{\vec{\boldsymbol{r}}}_{dma} \cdot {}^a\dot{\vec{\boldsymbol{r}}}_{dma} dm \\
\end{equation}where $S$ denotes the surface of the buoy, Substituting Eq. (\ref{eq:31}) in Eq. (\ref{eq:32}) to get,
\begin{align}
 \mathcal{T}  = \frac{1}{2} \int_S  & \begin{bmatrix} C_{sa}  & -[C_{es} ({}^e\vec{\boldsymbol{r}}_{cs} + {}^e\vec{\boldsymbol{r}}_{dmc})]^\times  & C_{se}
\end{bmatrix} \dot{\vec{\mathbf{x}}} \nonumber \\ 
\cdot & \begin{bmatrix}
C_{sa}  & -[C_{es} ({}^e\vec{\boldsymbol{r}}_{cs} + {}^e\vec{\boldsymbol{r}}_{dmc})]^\times  & C_{se} 
\end{bmatrix} \dot{\vec{\mathbf{x}}} \; dm \label{eq:33} 
\end{align} Let $\mathcal{H} = [C_{es} ({}^e\vec{\boldsymbol{r}}_{cs} + {}^e\vec{\boldsymbol{r}}_{dmc})]^\times$ such that Eq. (\ref{eq:33}) can be reduced to:  \begin{align}
\mathcal{T} &= \frac{1}{2} \int_S   \dot{\vec{\mathbf{x}}}^T   \begin{bmatrix}
C_{sa}^T  \\ -\mathcal{H}^T  \\ C_{se}^T
\end{bmatrix}\begin{bmatrix}
C_{sa}  & -\mathcal{H}  & C_{se} 
\end{bmatrix}
 \dot{\vec{\mathbf{x}}} dm \label{eq:35} \\
   &= \frac{1}{2}  \int_S  {\begin{bmatrix}
{}^a\dot{\vec{\boldsymbol{r}}}_{sa}^T \\
\vec{\omega}_{as}^T \\
{}^e\dot{\vec{\boldsymbol{r}}}_{dmc}^T
\end{bmatrix}}^T 
\begin{bmatrix}
 \boldsymbol{1} & - C_{sa}^T \mathcal{H} &   C_{sa}^T C_{se} \\
 - C_{sa}^T \mathcal{H} &  \mathcal{H}\mathcal{H}^T & - \mathcal{H}^T  C_{se}  \\
  C_{sa}^T C_{se}  &  - \mathcal{H}^T  C_{se}  & \boldsymbol{1} 
 \end{bmatrix} \nonumber \\ &\qquad \qquad \times \begin{bmatrix}
{}^a\dot{\vec{\boldsymbol{r}}}_{sa} \\
\vec{\omega}_{as} \\
{}^e\dot{\vec{\boldsymbol{r}}}_{dmc}
\end{bmatrix} dm 
\end{align}

 Finally, the kinetic energy of the system is expressed as:
 \begin{align}
 \mathcal{T}   = \frac{1}{2}  \int_S & \lbrace 
{}^a\dot{\vec{\boldsymbol{r}}}_{sa}^T( {}^a\dot{\vec{\boldsymbol{r}}}_{sa} - C_{sa}^T \mathcal{H} \vec{\omega}_{as} + C_{sa}^T C_{se} {}^e\dot{\vec{\boldsymbol{r}}}_{dmc}) \nonumber \\ &+
 \vec{\omega}_{as}^T( - C_{sa}^T \mathcal{H} {}^a\dot{\vec{\boldsymbol{r}}}_{sa} +  \mathcal{H}\mathcal{H}^T \vec{\omega}_{as} - \mathcal{H}^T  C_{se} {}^e\dot{\vec{\boldsymbol{r}}}_{dmc} ) \nonumber \\
 & + {}^e\dot{\vec{\boldsymbol{r}}}_{dmc}^T ( C_{sa}^T C_{se} {}^a\dot{\vec{\boldsymbol{r}}}_{sa}  - \mathcal{H}^T  C_{se} \vec{\omega}_{as} +  {}^e\dot{\vec{\boldsymbol{r}}}_{dmc})   \rbrace dm \label{eq:19_5}
\end{align}

It can be proved that kinetic energy components associated with the multiplication of any two of ${}^a\dot{\vec{\boldsymbol{r}}}_{sa}$, $\vec{\omega}_{as}$ and ${}^e\dot{\vec{\boldsymbol{r}}}_{dmc}$ result in zero terms, then Eq. \eqref{eq:19_5} is reduced to:
\begin{align}
\mathcal{T}     &= \frac{1}{2}  \int_S \left\lbrace
{}^a\dot{\vec{\boldsymbol{r}}}_{sa}^T  {}^a\dot{\vec{\boldsymbol{r}}}_{sa} +
\vec{\omega}_{as}^T \mathcal{H}\mathcal{H}^T \vec{\omega}_{as}  +
 {}^e\dot{\vec{\boldsymbol{r}}}_{dmc}^T  {}^e\dot{\vec{\boldsymbol{r}}}_{dmc}
   \right\rbrace dm \nonumber \\ 
&=    \underbrace{ \frac{1}{2}  {}^a\dot{\vec{\boldsymbol{r}}}_{sa}^T \boldsymbol{m} {}^a\dot{\vec{\boldsymbol{r}}}_{sa}}_{\mathcal{T}_{x}}  + \underbrace{ \frac{1}{2} 
\vec{\omega}_{as}^T \boldsymbol{\mathcal{J}}_s \vec{\omega}_{as}}_{\mathcal{T}_{\omega}} +
\underbrace{ \frac{1}{2}\int_s  {}^e\dot{\vec{\boldsymbol{r}}}_{dmc}^T   {}^e\dot{\vec{\boldsymbol{r}}}_{dmc}  dm}_{\mathcal{T}_{\boldsymbol{s}}} 
\label{eq:47_4}
\end{align}
where $\boldsymbol{m} = \boldsymbol{1}m$ is the VSB mass matrix and $\boldsymbol{\mathcal{J}}_s = \boldsymbol{1}\mathcal{J}_s $ is the second moment of inertia. Considering Eq. (\ref{eq:47_4}), the following is noted:
\begin{enumerate}
\item The first term yields the translational kinetic energy ${\mathcal{T}_{x}}$ of the rigid buoy.
\item The second term describes the rotational kinetic energy $\mathcal{T}_\omega$ of the rigid shell.
\item The third term accounts for the kinetic energy associated with the deformation of the shell ${\mathcal{T}_{\boldsymbol{s}}}$ and is simplified as follows:
\begin{align} 
\mathcal{T}_{\boldsymbol{s}} &= \frac{1}{2} \int_s {}^e\dot{\vec{\boldsymbol{r}}}_{dmc}^T {}^e\dot{\vec{\boldsymbol{r}}}_{dmc} dm \nonumber \\
&= \frac{1}{2}  \rho h \int_0^\pi \int_0^{2 \pi} {}^e\dot{\vec{\boldsymbol{r}}}_{dmc}^T {}^e\dot{\vec{\boldsymbol{r}}}_{dmc}r^2 \sin{\phi}  d\theta   d\phi  \label{eq:20_5} \\
&=\frac{1}{2} \rho h \int_0^\pi \int_0^{2 \pi} \begin{bmatrix} \dot{u} & 0 & \dot{v} \end{bmatrix}   \begin{bmatrix} \dot{u} \\ 0 \\ \dot{v} \end{bmatrix}  r^2 \sin{\phi}  d\theta   d\phi \nonumber \\  
  &= \frac{1}{2} (2 \pi \rho h) \int_0^\pi ( \left\lbrace \dot{u}^2 + \dot{v}^2   \right\rbrace  r^2 \sin{\phi} )  d\phi \label{eq:22_5}
 \end{align} 
\end{enumerate}
where $h$ and $\rho$ are the shell's thickness and material density, respectively. Eqs. \eqref{eq:22_5} describes the kinetic energy for spherical axisymmetric homogeneous thin shells. 


\subsection{Potential Energy for Flexible Spherical Buoys} \label{subsec:PE}
Considering the buoy as a spherical shell, the gravitational energy can be expressed as:
\begin{equation}
\mathcal{G}  = \int_s -{\vec{\boldsymbol{g}}} \cdot \vec{\boldsymbol{r}}_{sa} dm = mg  {\hat{\boldsymbol{a}}_3}^T \vec{\boldsymbol{r}}_{sa}   \label{eq:29}
\end{equation}
where $\vec{\boldsymbol{g}} = -g \hat{\boldsymbol{a}}_3$.
The strain energy-displacement expressions (membrane strains) for axisymmetric shells can be written as \cite{rao2007vibration,soedel2005vibrations,robert2016modeling,hogan2015dynamic}:
\begin{align}
\varepsilon_{\phi\phi} & = \frac{1}{r} \left(  \frac{\partial u}{\partial\phi}+v  \right) \label{eq:15.1} \\ 
\varepsilon_{\theta\theta} & = \frac{1}{r} \left( u \cot(\phi) + v \right) \label{eq:16.1}
\end{align}
and the total strain energy can be expressed as \cite{rao2007vibration,shabara2022bang}:
\begin{align}
\mathcal{U}_{\boldsymbol{s}}  = \frac{1}{2} \frac{Eh}{1-\nu^2} \int_0^{2 \pi} \int_0^{\pi} &\left\lbrace \varepsilon_{\phi \phi}^2+ \varepsilon_{\theta \theta}^2 + 2 \nu \varepsilon_{\phi\phi} \varepsilon_{\theta \theta} \right\rbrace \nonumber \\ &\times  r^2 \sin(\phi) d\phi \label{eq:18} 
\end{align}

Combining the sphere elastic and gravitational potential energies in Eqs.~\eqref{eq:29} and \eqref{eq:18} yields the total potential energy of the spherical shell buoy:
\begin{align}
\mathcal{\pi} &=  \mathcal{U}  + \mathcal{G} \nonumber \\ &=  \frac{1}{2} \frac{Eh}{1-\nu^2} \int_s \left\lbrace{ \varepsilon_e^{\phi \phi}}^2+ {\varepsilon_e^{\theta \theta}}^2 + 2 \nu {\varepsilon_e^{\phi\phi}} \varepsilon_e^{\theta \theta} \right\rbrace  r^2 \sin(\phi) dm \nonumber \\  & \qquad \qquad \qquad \qquad  + mg  \boldsymbol{1}_3^T \vec{\boldsymbol{r}}_{sa} \label{eq:45_1}   \\
&=  \frac{1}{2} \frac{Eh}{1-\nu^2} \int_0^{2 \pi} \int_0^\pi \left\lbrace{ \varepsilon_e^{\phi \phi}}^2+ {\varepsilon_e^{\theta \theta}}^2 + 2 \nu {\varepsilon_e^{\phi\phi}} \varepsilon_e^{\theta \theta} \right\rbrace \nonumber \\ & \qquad \qquad \qquad \qquad \times  r^2 \sin(\phi) d \theta d \phi + mg  \boldsymbol{1}_3^T \vec{\boldsymbol{r}}_{sa}
\end{align}

\section{Unconstrained Equations of Motion For Flexible Spherical Buoys} \label{sec:EOM}

The unconstrained equations of motion are here derived as a first step towards writing the constrained equations of motion. The Lagrangian for this buoy system can be written as the summation of three quantities \cite{junkins2009analytical}:
\begin{align}
    \mathcal{L} =\mathcal{L}_D + \mathcal{L}_B + \int_s \hat{\mathcal{L}} d \phi
\end{align}
where $\mathcal{L}_D(t,\vec{\mathbf{x}},\dot{\vec{\mathbf{x}}})$ and $\mathcal{L}_B$ are related to the discrete coordinates and the boundaries, respectively. , $\hat{\mathcal{L}}$ is the Lagrangian density function and it is a function of the discrete and distributed parameter coordinates. In this work there are no boundary terms in the Lagrangian equation, i.e., $\mathcal{L}_B = \boldsymbol{0}$. First, we will derive the equations of motion related to the discrete coordinates, then the equations of motion related to the distributed parameter coordinates will be derived where Rayleigh-Ritz approximation will be summoned.

\subsection{Equation of Motion Associated with Discrete Coordinates}
The Lagrangian for the discrete coordinates is expressed as:
\begin{align} 
\mathcal{L}_D = \mathcal{T}_D - \mathcal{\pi}_ D
            =  \frac{1}{2}   {}^a\dot{\vec{\boldsymbol{r}}}_{sa}^T \boldsymbol{m} {}^a\dot{\vec{\boldsymbol{r}}}_{sa} + \frac{1}{2}  {\vec{\boldsymbol{\omega}}}_{as}^T \boldsymbol{\mathcal{J}}_s {\vec{\boldsymbol{\omega}}}_{as} -mg  \boldsymbol{1}_3^T \vec{\boldsymbol{r}}_{sa} \label{eq:30_5}
\end{align}

The Lagrange Equation for the discrete coordinates is expressed as:
\begin{equation} \label{eq:49_3}
\frac{d}{dt}\left(\frac{\partial \mathcal{L}_D}{\partial \dot{\vec{\mathbf{x}}}}\right) - \frac{\partial \mathcal{L}_D}{\partial{\vec{\mathbf{x}}}}= \boldsymbol{0}
\end{equation}

To write the equations of motion of the discrete coordinates, we first write:
\begin{align}
\frac{\partial \mathcal{L}_D}{\partial{\vec{\boldsymbol{r}}_{sa}}} &= - mg\boldsymbol{1}_3 \\
\frac{d}{dt} \left( \frac{\partial \mathcal{L}_D }{\partial \dot{\vec{\boldsymbol{r}}}_{sa}} \right)   &=  \boldsymbol{m} {}^a\ddot{\vec{\boldsymbol{r}}}_{sa} \\
\frac{d}{dt}  \left( \frac{\partial \mathcal{L}_D }{\partial {\vec{\boldsymbol{\omega}}}_{as}}\right)
&=  \boldsymbol{\mathcal{J}}_s \dot{\vec{\boldsymbol{\omega}}}_{as} + \dot{\boldsymbol{\mathcal{J}}}_s {\vec{\boldsymbol{\omega}}}_{as} 
\end{align} 
The equations of motion for the translation and rotational motions become:
\begin{align}
     \boldsymbol{m} {}^a\ddot{\vec{\boldsymbol{r}}}_{sa} + mg\boldsymbol{1}_3 = 0 \label{eq:34_5} \\
      \boldsymbol{\mathcal{J}}_s \dot{\vec{\boldsymbol{\omega}}}_{as} + \dot{\boldsymbol{\mathcal{J}}}_s {\vec{\boldsymbol{\omega}}}_{as} = 0 \label{eq:35_5}
\end{align} The equations of motion described by Eq. \eqref{eq:34_5} and \eqref{eq:35_5} can be extended to include damping coefficient matrices as follows:

\begin{align}
     \boldsymbol{m} {}^a\ddot{\vec{\boldsymbol{r}}}_{sa} + \boldsymbol{D}_x {}^a\dot{\vec{\boldsymbol{r}}}_{sa} + mg\boldsymbol{1}_3 = \boldsymbol{0} \label{eq:36_5}\\
      \boldsymbol{\mathcal{J}}_s \dot{\vec{\boldsymbol{\omega}}}_{as} + (\dot{\boldsymbol{\mathcal{J}}}_s + \boldsymbol{D}_\omega ) {\vec{\boldsymbol{\omega}}}_{as}  = \boldsymbol{0} \label{eq:37_5}
\end{align}

where $\boldsymbol{D}_x$ and $\boldsymbol{D}_\omega$ are the damping matrices for the translation and rotation motions, respectively. 

\subsection{Equations of Motion Associated with The Distributed Parameter Coordinates}
The distributed parameters Lagrangian, $\hat{\mathcal{L}}$, is expressed as:
\begin{align}
    \hat{\mathcal{L}} &=   \frac{1}{2}\lbrace  (2 \pi \rho h)  \left\lbrace \dot{u}^2 + \dot{v}^2  \right\rbrace  r^2 \sin{\phi}   \rbrace \nonumber \\ 
    &\;\;\;\;\;\;-  \frac{1}{2} \frac{2 \pi Eh}{(1-\nu^2)} \Biggl\{ { \left(  \frac{\partial u}{\partial\phi}+v  \right) }^2+ {\left( u \cot(\phi) + v\right) }^2 \nonumber \\ 
    & \;\;\;\;\;\;+ 2 \nu {\left(  \frac{\partial u}{\partial\phi}+v  \right)} \left( u \cot(\phi) + v\right)  \Biggl\} r^2 \sin{\phi}
\end{align} and, the Lagrange equation for the distributed parameters is \cite{junkins2009analytical}: \begin{align}
    \frac{d}{dt} \left( \frac{\partial \hat{\mathcal{L}}}{\partial \dot{\vec{\boldsymbol{r}}}_{dmc}} \right) - \frac{\partial \mathcal{L}}{\partial \vec{\boldsymbol{r}}_{dmc}} &+ \frac{\partial }{\partial \vec{\boldsymbol{x}}} \left(\frac{\partial \hat{\mathcal{L}}}{\partial {\vec{\boldsymbol{r}}}_{dmc}'} \right) \nonumber \\ &\qquad  +  \frac{\partial^2 }{\partial \vec{\boldsymbol{x}}^2 } \left(\frac{\partial \hat{\mathcal{L}}}{\partial {\vec{\boldsymbol{r}}}_{dmc}''} \right)  = \boldsymbol{0}
\end{align} where $\vec{\boldsymbol x}$ is the vector of spatial coordinates in the $\hat{\boldsymbol{a}}$ directions, such that the transformation from the cartesian coordinates to spherical coordinates is as follows: \begin{align}
x_1 = r \sin{\phi} \leftrightarrows \frac{\partial}{\partial x_1} = \frac{1}{r \cos{\phi} } \frac{\partial}{\partial \phi} \\ 
x_3 = r \cos{\phi} \leftrightarrows  \frac{\partial}{\partial x_3} = \frac{-1}{r \sin{\phi} } \frac{\partial}{\partial \phi} 
\end{align} For ${{\vec{r}}_{dmc}}_1= u$: \begin{align}
&\frac{d}{dt} \left( \frac{\partial \hat{\mathcal{L}}}{\partial {\dot{u}}} \right) = 2 \pi \rho h r^2 \sin{\phi}  \ddot{u} \\
 &   \frac{\partial \hat{\mathcal{L}}}{\partial u} = -  \frac{2 \pi Eh}{(1-\nu^2)}  \cot(\phi) \Biggl\{ {  \left( u \cot(\phi) + v\right) } +  \nu   {\left(  \frac{\partial u}{\partial\phi}+v  \right)}   \Biggl\}\nonumber \\  &\qquad \qquad \times r^2 \sin{\phi} \\
& \frac{1}{r \cos{\phi} }  \frac{\partial }{\partial \phi }  \left(\frac{\partial \hat{\mathcal{L}}}{\partial {u}'} \right) = -   \frac{ r \pi Eh}{(1-\nu^2)}  \Biggl\{ {\left(  \frac{\partial^2 u}{\partial\phi^2}+ \frac{\partial v}{\partial\phi}  \right) }    \nonumber \\ & \qquad \qquad+   \nu \left(  \frac{\partial u}{\partial\phi}   \cot(\phi) + u \csc(\phi)  + v\right)  \Biggl\}  - \frac{ \pi Eh}{(1-\nu^2)} \nonumber\\& \qquad \qquad \times
    \Biggl\{ \left(  \frac{\partial u}{\partial\phi} +v  \right) + \nu \left( u \cot(\phi) + v\right)  \Biggl\} r \tan{\phi}
\end{align}

For ${{\vec{r}}_{dmc}}_3= v$: \begin{align}
 &  \frac{d}{dt} \left( \frac{\partial \hat{\mathcal{L}}}{\partial {\dot{v}}} \right) =  2 \pi \rho h r^2 \sin{\phi}  \ddot{v}  \\
&\frac{\partial \hat{\mathcal{L}}}{\partial v}  = -  \frac{2 \pi Eh}{(1-\nu^2)}  \Biggl\{ { \left(  \frac{\partial u}{\partial\phi}+v  \right) } + {\left(  u \cot(\phi) + v\right) } \nonumber \\ & \qquad+  \left( \nu {\left(  \frac{\partial u}{\partial\phi}+v  \right)} + \left( u \cot(\phi) + v \right) \right)  \Biggl\}  r^2 \sin{\phi} \\
&\frac{-1}{r \sin{\phi} } \frac{\partial}{\partial \phi}  \left(\frac{\partial \hat{\mathcal{L}}}{\partial {v}'} \right) = \frac{-1}{r \sin{\phi} } \frac{\partial}{\partial \phi}   \left( 0 \right) =  0
\end{align} Accordingly, the equations of motion associated with the distributed parameter are expressed as:

\begin{align}
 & 2 \pi  \rho h r^2 \sin{\phi}  \ddot{u}   - \frac{2 \pi Eh}{(1-\nu^2)}  \cot(\phi)   \Biggl\{ {  \left( u \cot(\phi) + v\right) } \nonumber  \\
& +  \nu   {\left(  \frac{\partial u}{\partial\phi} +v  \right)} \Biggl\}   r^2 \sin{\phi} -   \frac{ \pi Eh}{(1-\nu^2)} \times \nonumber  \\
&\left\lbrace{\left(  \frac{\partial^2 u}{\partial\phi^2}+ \frac{\partial v}{\partial\phi}  \right) }+   \nu \left(  \frac{\partial u}{\partial\phi} \cot(\phi) + u \csc(\phi)  + v\right)  \right\rbrace r \nonumber  \\
& -   \frac{r \pi Eh}{(1-\nu^2)} \left\lbrace{\left(  \frac{\partial u}{\partial\phi} +v  \right) }+   \nu \left( u \cot(\phi) + v\right)  \right\rbrace  \tan{\phi}= 0 \label{eq:49_6} \\
  2 \pi & \rho h r^2 \sin{\phi}  \ddot{v}   -  \frac{2 \pi Eh}{(1-\nu^2)}  \Biggl\{ {\left(  \frac{\partial u}{\partial\phi}+v  \right) }+ {\left(  u \cot(\phi) + v\right) } \nonumber \\ & + \left( \nu {\left(  \frac{\partial u}{\partial\phi}+v  \right)} + \left( u \cot(\phi) + v \right) \right)  \Biggl\} r^2 \sin{\phi}   = 0 \label{eq:50_5}
\end{align}

Equations \eqref{eq:49_6} and \eqref{eq:50_5} are the equations of motion of deformation of axisymmetric homogeneous spherical shell. These equations of motion can not be solved analytically, and an approximate method is used. Most of the equations of motions for continuous systems are usually difficult to obtain; this difficulty arises either from the difficulty in solving the governing equations or from imposing the boundary conditions. Here, an approximate method is implemented to convert the partial differential equations to ordinary differential equations.

\subsubsection*{Rayleigh-Ritz Approximation}
The Rayleigh-Ritz approximation method is applied in this article.
Each component of the displacement vector $\vec{\boldsymbol{r}}_{dmc}$ is assumed to have the following form \cite{robert2016modeling,hogan2015dynamic}::
\begin{align}
u(\phi,t) &= \sum_{n=1}^N \Psi_{n}^{\phi}(\phi)\eta_n(t) = \underbrace{[\Psi_{1}^{\phi} \ldots \Psi_{N}^{\phi }]}_{\boldsymbol{\Psi}_{e}^{\phi}}  \underbrace{\begin{bmatrix}
    \eta_{1}(t) \\
    \vdots \\
    \eta_{N}(t) 
\end{bmatrix}}_{\boldsymbol{\eta}(t)} \nonumber \\  &= \boldsymbol{\Psi}_{e}^{\phi}(\phi) \boldsymbol{\eta}(t) \label{eq:51_5} \\
v(\phi,t) &= \sum_{n=1}^N \Psi_{n}^{r}(\phi)\eta_n(t)= \underbrace{[\Psi_{1}^{r} \ldots \Psi_{n}^{r}]}_{\boldsymbol{\Psi}_{e}^{r}} 
\underbrace{\begin{bmatrix}
    \eta_{1}(t) \\
    \vdots \\
    \eta_{N}(t) 
\end{bmatrix}}_{\boldsymbol{\eta}(t)} \nonumber \\ &= \boldsymbol{\Psi}_{e}^{r}(\phi) \boldsymbol{\eta}(t) \label{eq:52_5}
\end{align}where the functions $\Psi_{n}^{\phi}$ and $\Psi_{n}^{r}$ are trial (admissible) functions of $\phi$ and the functions $\eta_n$ are functions of time $t$, $\forall$ $n=1, \cdots, N$.
Therefore, the displacement vector can be expressed in the $\hat{\boldsymbol{e}}$ frame as follows:
\begin{equation}\label{eq:5}
\vec{\boldsymbol{r}}_{dmc}(\phi,t) = \underbrace{\begin{bmatrix}
\boldsymbol{\Psi}_{e}^{\phi}(\phi) \\
\boldsymbol{0} \\
\boldsymbol{\Psi}_{e}^{r}(\phi) 
\end{bmatrix}}_{\boldsymbol{\Phi}_e} \boldsymbol{\eta}(t) = \boldsymbol{\Phi}_e(\phi) \boldsymbol{\eta}(t) \end{equation} The Legendre functions of the first kind $P_n$ \cite{kreyszig2009advanced} can serve as shape functions for the Ritz-Rayleigh method to satisfy the essential geometrical (Dirichlet) boundary conditions  \cite{naghdi1962vibrations,hogan2015dynamic,robert2016modeling,soedel2005vibrations,newman1994wave,raouf1990non}, as follows: \begin{align}
\Psi_{n}^{\phi}(\phi) &= A \frac{d P_n (\cos(\phi))}{d \phi} \text{, and} \\ \nonumber \Psi_{n}^{r}(\phi) &= A \frac{(1+(1+\nu))\Omega_n^2}{1-\Omega_n^2}   P_n (\cos(\phi))
\end{align} where the coefficients of the equations above form an eigenvector for the Legendre differential equation, i.e. the constant "A" can take any real value. $\Omega_n^2$ is a dimensionless frequency parameter expressed as \cite{soedel2005vibrations}:
\begin{equation} \label{eq:7}
\Omega^2_{n} = \frac{1}{2(1-\nu^2)} (A \pm \sqrt{A^2 -4mB})
\end{equation}where $\nu$ is the Poisson's ratio, and
\begin{align}
m &= n(n+1)-2,  \: n \in \mathbb{Z}^{+} \label{eq:10}\\
B &= 1+ \nu^2 +\frac{1}{12} [(m+1)^2 - \nu^2] \label{eq:9}\\
A &= 3(1+\nu) + m+\frac{1}{2} \left[ \frac{h}{r} \right]^2 (m+3)(m+1+\nu) \label{eq:8}
\end{align}
From \cite{soedel2005vibrations,naghdi1962vibrations} the natural frequencies in radians per second for spherical shells are calculated using Eq. (\ref{eq:11})
\begin{equation} \label{eq:11}
\omega_n^2 = \frac{E}{r^2 \rho \Omega_n^2}
\end{equation} where $E$ is Young's Modulus of elasticity, $r$ is the non-deformed radius of the shell, and $\rho$ is the density of the shell material.
When $n = 0$, the vibration mode corresponds to the breathing mode (volumetric or pulsating modes) which is a pure radial vibration mode \cite{soedel2005vibrations,nayfeh2006axisymmetric,raouf1990non}. For $n > 0$, the $\pm$ sign in Eq.~\eqref{eq:7} yields the modes corresponding to the membrane vibration modes and bending vibration modes (this is demonstrated in Fig. \ref{fig:modes} in the results section). The bending vibration modes are obtained when using the negative sign; these modes are sensitive to the $h/r$ ratio. On the other hand, the membrane modes are insensitive to the change in the $h/r$ ratio. To obtain the approximated equations of motion using the Rayleigh-Ritz method, the approximated displacement vector needs to be substituted in the kinetic and strain energy equations as follows.
The kinetic energy $\mathcal{T}_{\boldsymbol{s}}$ is approximated by substituting Eq.s \eqref{eq:51_5} and \eqref{eq:52_5} into Eq. \eqref{eq:22_5}; to get: 
\begin{align} 
\mathcal{T}_{\boldsymbol{s}} &= \frac{1}{2} \dot{\boldsymbol{\eta}}^T \underbrace{\Bigg\{ 2 \pi \rho h \int^{\pi}_{0} \left( \boldsymbol{\Psi}_{e}^{{\phi}^T} \boldsymbol{\Psi}_{e}^{\phi} +  \boldsymbol{\Psi}_{e}^{{r}^T} \boldsymbol{\Psi}_{e}^{r}\right)  r^2 \sin \phi d \phi \Bigg\}}_{\boldsymbol{M}_{ee}} \dot{\boldsymbol{\eta}} \label{eq:13} \\
& = \frac{1}{2} \dot{\boldsymbol{\eta}}^{T} \boldsymbol{M}_{ee} \dot{\boldsymbol{\eta}}  \label{eq:14}
\end{align}
The strain energy is approximated by first substituting Eqs. \eqref{eq:51_5} and \eqref{eq:52_5} in Eqs. \eqref{eq:15.1} and \eqref{eq:16.1} to get:
\begin{align}
\varepsilon_{\phi\phi} & =  \frac{1}{r} \left(  \frac{\partial \Psi_e^\phi}{\partial\phi}+\Psi_e^r\right)\boldsymbol{\eta} \label{eq:15} \\ 
\varepsilon_{\theta\theta} & =   \frac{1}{r} \left( \Psi_e^\phi \cot(\phi) + \Psi^r_e \right) \boldsymbol{\eta} \label{eq:16}
\end{align}
Then the strain energy equation in Eq. \eqref{eq:18} becomes:
\begin{align}
\mathcal{U}_{\boldsymbol{s}} & = \frac{1}{2}{\boldsymbol{\eta}}^T \Bigg[ \frac{2\pi Eh}{1-\nu^2} \int_0^{\pi} \Bigg\{ \left(    \frac{\partial \Psi_e^\phi}{\partial\phi}+\Psi_e^r\right)^T \left(  \frac{\partial \Psi_e^\phi}{\partial\phi}+\Psi_e^r\right) \nonumber \\
&+  \left( \Psi_e^\phi \cot(\phi) + \Psi^r_e \right)^T  \left( \Psi_e^\phi \cot(\phi) + \Psi^r_e \right) \nonumber \\
&+\nu \Biggl(  \left(  \frac{\partial \Psi_e^\phi}{\partial\phi}+\Psi_e^r\right)^T   \left( \Psi_e^\phi \cot(\phi) + \Psi^r_e \right) \nonumber \\
&+  \left( \Psi_e^\phi \cot(\phi) + \Psi^r_e \right)^T    \left(  \frac{\partial \Psi_e^\phi}{\partial\phi} + \Psi_e^r\right) \Biggl) \Bigg\}\sin(\phi)  d\phi \Bigg] \boldsymbol{\eta} \\
&= \frac{1}{2} {\boldsymbol{\eta}}^T \boldsymbol{K}_{ee} \boldsymbol{\eta} \label{eq:19}
\end{align}
where 

\begin{align}
\boldsymbol{K}_{ee}   &= \frac{2\pi Eh}{1-\nu^2} \int_0^{\pi}  \Biggl\{ \zeta_{\phi\phi}^T \zeta_{\phi\phi} + {\zeta_{\theta \theta}}^T \zeta_{\theta \theta} \nonumber  \\ 
&\qquad \qquad +\nu \left( \zeta_{\phi\phi}^T \zeta_{\phi\phi} +{\zeta_{\theta \theta}}^T \zeta_{\theta \theta}\right) \Biggl\}  \sin(\phi) d\phi \\
\zeta_{\phi\phi} &=   \left(  \frac{\partial \Psi_e^\phi}{\partial\phi}+\Psi_e^r\right) \label{eq:20}\\ 
\zeta_{\theta \theta} &= \left( \Psi_e^\phi \cot(\phi) + \Psi^r_e \right) \label{eq:21}
\end{align}

Using the above approximate expressions for the kinetic and strain energies, the Lagrangian $\mathcal{L}_s$ is expressed as follows:
\begin{align}
   \mathcal{L}_s & \equiv  \mathcal{T}_{\boldsymbol{s}}  - \mathcal{U}_{\boldsymbol{s}}   =
   \frac{1}{2} \dot{\boldsymbol{\eta}}^T \boldsymbol{M}_{ee} \dot{\boldsymbol{\eta}} - \frac{1}{2} {\boldsymbol{\eta}}^T \boldsymbol{K}_{ee} \boldsymbol{\eta} \label{eq:24_1}
\end{align}
With this approximation, the Lagrangian can be used to write the equations of motion of the flexible buoy in the form:
\begin{align}
   \frac{d}{dt}\left(\frac{\partial \mathcal{L}_s}{\partial\dot{\boldsymbol{\eta}}}\right) &-\frac{\partial \mathcal{L}_s}{\partial{\boldsymbol{\eta}}}= \boldsymbol{0} \label{eq:25_1}
\end{align}
From Eq.~\eqref{eq:24_1}, one can write:
\begin{align}
      \frac{\partial \mathcal{L}_s}{\partial{\boldsymbol{\eta}}} & =- \boldsymbol{K}_{ee} \boldsymbol{\eta}, \; \text{and} \; \frac{\partial \mathcal{L}_s}{\partial\dot{\boldsymbol{\eta}}}  = \boldsymbol{M}_{ee} \dot{\boldsymbol{\eta}}  \nonumber\\
\therefore   \frac{d}{dt}\left(\frac{\partial \mathcal{L}_s}{\partial\dot{\boldsymbol{\eta}}}\right) &=  \boldsymbol{M}_{ee} \ddot{\boldsymbol{\eta}}
\end{align}
Hence, the equation of motion for the shell of a flexible buoy vibrating axisymmetricly is:
\begin{align}
   \boldsymbol{M}_{ee} \ddot{\boldsymbol{\eta}} &+ \boldsymbol{K}_{ee} \boldsymbol{\eta} = \boldsymbol{0} \label{eq:26_1}
\end{align}

Equation (\ref{eq:26_1}) can be further extended to include Rayleigh damping as follows:
\begin{equation} \label{eq:23}
\boldsymbol{M}_{ee} \ddot{\boldsymbol{\eta}} + \boldsymbol{D}_{ee} \dot{\boldsymbol{\eta}} + \boldsymbol{K}_{ee} \boldsymbol{\eta} = \boldsymbol{0} 
\end{equation}  
where $\boldsymbol{D}_{ee}$ is a {proportional} damping matrix which is a function of the mass and stiffness matrices as follows:
\begin{equation} \label{eq:23_1}
\boldsymbol{D}_{ee} = \alpha_d \boldsymbol{M}_{ee} + \beta_d \boldsymbol{K}_{ee} \nonumber
\end{equation}  
where the $\alpha_d$ and $\beta_d$ are real scalars named the mass and stiffness matrix multipliers with units 1/sec and sec, respectively \cite{alipour2008study,liu1995formulation,thompson2017ansys}. 
Combining the equations of motion from Eqs. \eqref{eq:34_5}, \eqref{eq:35_5} and  \eqref{eq:23} yields the equation of motion of the flexible buoy:
\begin{align}
& \boldsymbol M {\ddot{\vec{\mathbf{x}}}}  + ( \dot{\boldsymbol{M}} + \boldsymbol{D}) {\dot{\vec{\mathbf{x}}}}+ \begin{bmatrix}
mg \boldsymbol{1}_3^T & \boldsymbol{0} &  (\boldsymbol{K}_{ee} \boldsymbol{\eta})^T
\end{bmatrix}^T = \boldsymbol{0} \label{eq:50_2}
\end{align} 
where
\begin{align}
\dot{\vec{\mathbf{x}}} &= \begin{bmatrix}
{}^a\dot{\vec{\boldsymbol{r}}}_{sa}^T & {\vec{\mathbf{\omega}}}_{sa}^T  &
\dot{\boldsymbol{\eta}}^T  \end{bmatrix}^T_{(6+N) \times 1}     \\
\boldsymbol M  &=  \text{diag} \{ \begin{matrix} 
       \boldsymbol{m} & \boldsymbol{\mathcal{J}}_s & \boldsymbol{M}_{ee} \end{matrix} \}_{(6+N) \times (6+N)}  \\
\dot{\boldsymbol{M}}    &= \text{diag} \{ \begin{matrix} 
       \boldsymbol{0} & \dot{\boldsymbol{\mathcal{J}}}_s &  \boldsymbol{0}\end{matrix} \}_{(6+N) \times (6+N)}     \\
\boldsymbol{D} &= \text{diag} \{ \begin{matrix} 
     \boldsymbol{D}_x & \boldsymbol{D}_{\omega} & \boldsymbol{D}_{ee} \end{matrix} \}_{(6+N) \times (6+N)} 
\end{align}

It is noted here that the system mass matrix is a function of time because $\boldsymbol{\mathcal{J}}_s$ is a function of time. The vector $\vec{\boldsymbol{\omega}}_{sa}$ describes the instantaneous body angular velocities in the body frame $\hat{\boldsymbol{s}}$ with respect to the inertial frame $\hat{\boldsymbol{a}}$. To avoid integrating the direction cosine matrices, the Euler angles are used for orientation (attitude) descriptions. 
Let  $[\boldsymbol{B}(\theta)]$ be the mapping matrix that converts the angular velocity $\vec{\boldsymbol{\omega}}_{sa}$ to  Euler angle rates $\dot{\boldsymbol\theta}_{sa}$; hence we can write:
\begin{align} \small
\vec{\boldsymbol{\omega}}_{sa} &= [\boldsymbol{B}(\theta)]^{-1} \dot{\boldsymbol\theta}_{sa} \\
\dot{\vec{\boldsymbol{\omega}}}_{sa} &= [\boldsymbol{B}(\theta)]^{-1} \ddot{\boldsymbol\theta}_{sa} + [\dot{\boldsymbol{B}}(\theta)]^{-1} \dot{\boldsymbol\theta}_{sa} \label{eq:52_2} \\
[\boldsymbol{B}(\theta)]^{-1} &= \begin{bmatrix}
-\sin{\theta_2} & 0 & 1 \\
\cos{\theta_2} \sin{\theta_3} & \cos{\theta_3}& 0\\
\cos{\theta_2} \cos{\theta_3} & -\sin{\theta_3}  & 0 \\
\end{bmatrix} \label{eq:53_2} \end{align} \begin{align}
&[\dot{\boldsymbol{B}}(\theta)]^{-1}  \nonumber \\ &=- \begin{bmatrix}
 \dot{\theta}_2 \cos{\theta_2}& 0 & 0 \\
\dot{\theta}_2 \sin{\theta_2} \sin{\theta_3} - \dot{\theta}_3 \cos{\theta_2} \cos{\theta_3} &  \dot{\theta}_3\sin{\theta_3}& 0\\
\dot{\theta}_2 \sin{\theta_2} \cos{\theta_3} + \dot{\theta}_3
\cos{\theta}_2 \sin{\theta_3} & \dot{\theta}_3 \cos{\theta_3}  & 0 \\
\end{bmatrix} \label{eq:54_2}
\end{align}

where Eqs. \eqref{eq:53_2} and \eqref{eq:54_2} are derived using a 3-2-1 Euler angles sequence. Substituting Eq. \eqref{eq:52_2} into the inertia term in \eqref{eq:50_2} to get:

\begin{align}
\boldsymbol M & \begin{bmatrix}
{}^a\ddot{\vec{\boldsymbol{r}}}_{sa} \\
\dot{\vec{\mathbf{\omega}}}_{sa}  \\
\ddot{\boldsymbol{\eta}}  \end{bmatrix} + (\dot{\boldsymbol{M}}+\boldsymbol{D}  ) \begin{bmatrix}
{}^a\dot{\vec{\boldsymbol{r}}}_{sa} \\
\vec{\mathbf{\omega}}_{sa}  \\
\dot{\boldsymbol{\eta}}  \end{bmatrix}+ \boldsymbol{a}  \nonumber\\ 
= \boldsymbol M & \begin{bmatrix}
{}^a\ddot{\vec{\boldsymbol{r}}}_{sa} \\
[\boldsymbol{B}(\theta)]^{-1}  \ddot{\boldsymbol\theta}_{sa}  \\
\ddot{\boldsymbol{\eta}}^{ee} \end{bmatrix}
 + \boldsymbol M \begin{bmatrix} \boldsymbol{0}  \\
[\dot{\boldsymbol{B}}(\theta)]^{-1}  \dot{\boldsymbol\theta}_{sa} \\\boldsymbol{0} 
\end{bmatrix} +  \boldsymbol{D} \begin{bmatrix}
{}^a\dot{\vec{\boldsymbol{r}}}_{sa} \\
\vec{\mathbf{\omega}}_{sa}  \\
\dot{\boldsymbol{\eta}}  \end{bmatrix} \nonumber \\
& + \begin{bmatrix}
mg\boldsymbol{1}_3 \\ \dot{\boldsymbol{\mathcal{J}}}_s   \vec{\mathbf{\omega}}_{sa}\\ \boldsymbol{K}_{ee}\boldsymbol{\eta}
\end{bmatrix}\nonumber  \\
= \boldsymbol M & \bar{\boldsymbol{B}} \ddot{\mathbf{q}} + \boldsymbol M \dot{\bar{\boldsymbol{B}}} \dot{\mathbf{q}} + \boldsymbol{D} \dot{\mathbf{q}}  + \boldsymbol{a}
\end{align} where $\mathbf{q}$ is the generalized coordinates vector and is defined by Eq. \eqref{eq:57_2}:
\begin{align}
\mathbf{q} &= \begin{bmatrix} \vec{\boldsymbol r}_{sa}^T & \boldsymbol\theta_{sa}^T & \boldsymbol{\eta}^T  \end{bmatrix}^T \label{eq:57_2}
 \end{align} where
\begin{align}
\bar{\boldsymbol B} &= \text{diag} \{\begin{matrix} 
 \boldsymbol 1 &[\boldsymbol{B}(\theta)]^{-1} & \boldsymbol 1 \end{matrix}\}, \\
 \dot{\bar{\boldsymbol B}} &= \text{diag} \{\begin{matrix}
\boldsymbol{0} & [\dot{\boldsymbol{B}}(\theta)]^{-1} & \boldsymbol{0}
\end{matrix}\},\\
\boldsymbol{a} &= {\begin{bmatrix}
mg\boldsymbol{1}_3^T & ( \dot{\boldsymbol{\mathcal{J}}}_s \vec{\mathbf{\omega}}_{sa})^T & (\boldsymbol{K}_{ee}\boldsymbol{\eta})^T
\end{bmatrix}}^T 
 \end{align}
Then the Equation of motion for free buoy is expressed as: 
\begin{equation}
\bar{\boldsymbol{B}}^T \boldsymbol M \bar{\boldsymbol{B}}\ddot{\mathbf{q}} + \bar{\boldsymbol B}^T (\boldsymbol M \dot{\bar{\boldsymbol B}}\dot{\mathbf{q}} + \boldsymbol{D} {\bar{{\boldsymbol B}}} \dot{\mathbf{q}} + \mathbf{a}) = \boldsymbol{0} \label{eq:60_2}
\end{equation}
The $\mathbf{q}$ is a $(6+N)$ column vector of generalized coordinates, $\boldsymbol{D}$ is the system damping matrix, $\boldsymbol{D}_x$ and  $\boldsymbol{D}_{\omega}$ are the transnational and rotational damping coefficients matrices, respectively.

\section{Forced Constrained Equations of Motion} \label{sec:const_forced}

The Lagrange equation for the forced motion for the discrete coordinates is expressed as:
\begin{equation} \label{eq:49_5}
\frac{d}{dt}\left(\frac{\partial \mathcal{L}}{\partial \dot{\mathbf{q}}}\right) - \frac{\partial \mathcal{L}}{\partial{\mathbf{q}}}= \boldsymbol{Q}
\end{equation}
where, $\boldsymbol Q$ is a column vector of generalized forces. The external generalized forces on a buoy are the generalized $PTO$ force ``$\boldsymbol Q_{pto}$", generalized buoyant force ``$\boldsymbol Q_{b}$", generalized radiation forces $Q_r$, and generalized hydrodynamic excitation forces ``$\boldsymbol Q_e$".

The WECs can have an internal structure to install a $PTO$ unit at an angle $\phi_c$ on the shell. An example of this the VSB WEC tested in \cite{shabara2021numerical}, where a concert plate was attached at an angle $\phi_c = 90^o$ to apply the PTO force. 
The elastic surface displacement $\vec{\boldsymbol{r}}_{dmc}$ at $\phi_c = 90^o$ should equal to zero, resulting in the holonomic constraint in Eq. (\ref{eq:54}): \begin{equation} \label{eq:54}
{\vec{\boldsymbol{r}}}_{dmc} (\phi_c) = \boldsymbol{\Phi}_c \boldsymbol{\eta} = \boldsymbol{0}
\end{equation} where ${\boldsymbol{\Phi}}_c = {\boldsymbol{\Phi}}_{e}(\phi_c)$.

Hence, the most general form for the Equation of motion for a Spherical VSB is formed by equating the free vibration equation (Eq. \eqref{eq:60_2}) to the external and constraint generalized forces as shown in Eq. \eqref{eq:87_5}:
\begin{align}
\bar{\boldsymbol{B}}^T \boldsymbol M \bar{\boldsymbol{B}}\ddot{\mathbf{q}} + \bar{\boldsymbol B}^T (\boldsymbol M \dot{\bar{\boldsymbol B}}\dot{\mathbf{q}} &+ \boldsymbol{D} {\bar{{\boldsymbol B}}} \dot{\mathbf{q}} + \mathbf{a})  = \boldsymbol Q_c  + \boldsymbol Q \nonumber \\
&= \boldsymbol A^T \boldsymbol{\lambda} + \boldsymbol Q_{pto}+ \boldsymbol Q_{b} + \boldsymbol Q_e \label{eq:87_5}
\end{align}
where $\boldsymbol{Q}_c= \boldsymbol A^T \boldsymbol{\lambda}$ is the generalized constraint vector, $\boldsymbol \lambda$ column matrix of Lagrange multipliers with dimensions ($3\times1$), and $\boldsymbol A$ is the Jacobian constraint matrix.

 To write the constraint in Eq.~\eqref{eq:54} in the generalized coordinates, we write the following transformation:
\begin{align} \label{eq:55}
\boldsymbol{0} = {\boldsymbol{\Phi}}_c \boldsymbol{\eta} \equiv \underbrace{\begin{bmatrix}
   \boldsymbol{0} & \boldsymbol{0} & {\boldsymbol{\Phi}}_c \end{bmatrix}}_{\bar{\boldsymbol{\Phi}}_c } \bar{\boldsymbol{B}} \mathbf{q} = \underbrace{\bar{\boldsymbol{\Phi}}_c \bar{\boldsymbol{B}}}_{\boldsymbol{A}} \mathbf{q} 
\end{align}

Eq. \eqref{eq:55} can be written as $\boldsymbol{A} \mathbf{q} = \boldsymbol{0}$, where $\boldsymbol{A} = \bar{\boldsymbol{\Phi}}_c \bar{\boldsymbol{B}}$. 
Reference \cite{udwadia2000fundamental} shows that the generalized constraint forces for ideal constraints can be expressed as:
\begin{equation} \label{eq:61}
  \boldsymbol Q_c =  \boldsymbol M^{\frac{1}{2}}(\boldsymbol A \boldsymbol M^{-\frac{1}{2}})^\dagger (\boldsymbol b - \boldsymbol{AM}^{-1} (\boldsymbol Q-\bar{\boldsymbol a}))
\end{equation}
where $\boldsymbol b = - (\ddot{\boldsymbol A} \mathbf{q} +  2 \dot{\boldsymbol A} \dot{\mathbf{q}})$, and $\bar{\boldsymbol a} = \bar{\boldsymbol B }^T (\boldsymbol M \dot{\bar{\boldsymbol B}}\dot{\mathbf{q}} + \boldsymbol{D} {\bar{{\boldsymbol B}}} \dot{\mathbf{q}} + \mathbf{a})$. 

It is common in WEC analysis to assume that the WEC is only heaving (moving only in the vertical direction,) and hence we can simplify the equations by assuming that the buoy is not rotating; that is $\dot{\bar{\boldsymbol B}} = \boldsymbol 0$. In such case, the vector $\boldsymbol b = \boldsymbol 0$.
Next, the expression for each of the external forces is developed.

\subsection{Power Take-off Unit Force} 
The PTO can be either active or passive \cite{shabara2021numerical}, in this paper the PTO is assumed a passive damping. To apply the PTO effect for that case one can either include its effect in the $\boldsymbol{D}_x$ matrix in Eq. \eqref{eq:87_5}, or as an external force in Eq. \eqref{eq:73_3}.

\begin{equation} \label{eq:73_3}
    \vec{\mathbit{f}}^{pto} =   -c \: {}^a\dot{\vec{\mathit{\boldsymbol{r}}}}_{sa,3} \: \boldsymbol{a}_3
\end{equation}

Where $\vec{\mathbit{f}}^{pto}$ is the damping force and $c$ is the damping coefficient. Note that the equation above applies to heave only motion. To compute the generalized force, the general transformation takes the form \cite{junkins2009analytical}: 
\begin{equation} \label{eq:60_1g} 
  \mathbit{Q}_j =  \sum_{i=1}^{6+N} f_i \cdot \frac{\partial \vec{ r_i}}{\partial q_j}
\end{equation}
Hence, the generalized force corresponding to the PTO force takes the form:
\begin{equation} \label{eq:60_1} 
  \mathbit{Q}_j^{pto} =  \vec{\boldsymbol{f}}^{pto} \cdot \frac{\partial \vec{r}_{sa}}{\partial q_j}, \, \, \forall j=1, 2, \cdots 6+N
\end{equation}

We can then write the PTO generalized forces using Eq. \eqref{eq:57_2} and Eq. \eqref{eq:73_3} as follows:
\begin{align*}
\mathbit{Q}_1^{pto} &=  \vec{\boldsymbol{f}}^{pto}   \cdot \frac{\partial  \vec{\boldsymbol{r}}_{sa} }{\partial \vec{\boldsymbol{r}}_{sa,1}} = -c {}^a\dot{\vec{\mathit{\boldsymbol{r}}}}_{sa,3} \boldsymbol{a}_3 \cdot \boldsymbol{1}_1 =0  \\
 \mathbit{Q}_2^{pto} &=  \vec{\boldsymbol{f}}^{pto}   \cdot \frac{\partial  \vec{\boldsymbol{r}}_{sa} }{\partial \vec{\boldsymbol{r}}_{sa,2}} = -c {}^a\dot{\vec{\mathit{\boldsymbol{r}}}}_{sa,3}  \boldsymbol{a}_3 \cdot \boldsymbol{1}_2 =0 \\
 \mathbit{Q}_3^{pto} &= \vec{\boldsymbol{f}}^{pto}   \cdot \frac{\partial  \vec{\boldsymbol{r}}_{sa} }{\partial \vec{\boldsymbol{r}}_{sa,3}} = -c {}^a\dot{\vec{\mathit{\boldsymbol{r}}}}_{sa,3}  \boldsymbol{a}_3 \cdot \boldsymbol{1}_3 = -c {}^a\dot{\vec{\mathit{\boldsymbol{r}}}}_{sa,3} \nonumber \\
 \mathbit{Q}_4^{pto} &= \vec{\boldsymbol{f}}^{pto} \cdot \frac{\partial  \vec{\boldsymbol{r}}_{sa} }{\partial \theta_1 } = 0, \:\:\:\:
 \mathbit{Q}_5^{pto} = \vec{\boldsymbol{f}}^{pto} \cdot \frac{\partial  \vec{\boldsymbol{r}}_{sa} }{\partial \theta_2 } = 0,  \:\:\:\: \\
 \mathbit{Q}_6^{pto} &= \vec{\boldsymbol{f}}^{pto} \cdot \frac{\partial  \vec{\boldsymbol{r}}_{sa} }{\partial {\eta_1} } = 0, \:\:\:
\hdots \hdots \:\:\:
 \mathbit{Q}_{6+N}^{pto} =  \vec{\boldsymbol{f}}^{pto} \cdot \frac{\partial  \vec{\boldsymbol{r}}_{sa} }{\partial {\eta_N} } = 0
\end{align*} 
In a compact form, 
\begin{equation}
    \mathbit{Q}^{pto} = - \begin{bmatrix}
      c \: {}^a\dot{\vec{\mathit{\boldsymbol{r}}}}_{sa,3} \: \boldsymbol{1}_3 \\ \boldsymbol{0} \\ \boldsymbol{0}
    \end{bmatrix}
\end{equation}

\subsection{Hydro Forces} \label{subsec:3.1}
This work describes the method of applying the hydrodynamic loads on the shells in the Lagrangian mechanics realm, the hydrodynamic loads can be obtained from BEM solvers or CFD solvers (ex., capytaine, NEMOH and openFOAM).

In this work we assume that the average of the summation of hydrodynamic pressures (excitation, radiation, etc) is sinusoidal. It is assumed that the average hydrodynamic pressures $P_{hydrod}$ is uniformly distributed around the submerged buoy's volume; such that the hydrodynamic force on a submerged partition is expressed as ${f_{hydrod}}_i = P_{hydrod} {A}_i$, where $A_i$ is the $i^{th}$ partition surface area, the calculation of the surface area is detailed in section \ref{subsec:sub_area_calc}. 

The hydrostatic force (buoyant force) $\vec{\mathbit{f}}^b$ is exerted on the buoy's surface due to the displaced water by the submerged volume. For the case of a FSB, the buoyant force acts on the center of buoyancy along the $\hat{\boldsymbol{a}}_3$ direction. On the other hand, for a VSB, the buoyant force is computed as the integration of pressure over the surface. This buoyant pressure contributes to the shell's deformation; hence Eq. \eqref{eq:68} is used to compute the buoyant force at each node $i$:
\begin{eqnarray}
\vec{\mathbit{f}}^b_i =
    \frac{\rho_w V_{s,i} g }{ \cos (\pi - \psi_i)} {\hat{\boldsymbol{c}}_3^T}_i \label{eq:68} 
\end{eqnarray}where $\rho_w$ is the density of the water, $V_{s,i}$ is the submerged volume corresponding to the node $i$, and the angle $\phi$ is defined as shown in Fig.~\ref{fig:1}. The calculation of the buoyant force for the VSB takes into account the continuous change in the buoy shape and the submerged volume. 
In the rest of the current subsection, the methodology followed to calculate the submerged volume is presented, and then the generalized hydrodynamic force is calculated.

\subsubsection{Submerged Volume Calculation} \label{subsec:sub_volume_calc}
\begin{figure}%
    \centering
    \subfloat[\centering Submerged Nodes]{{\includegraphics[scale=0.22]{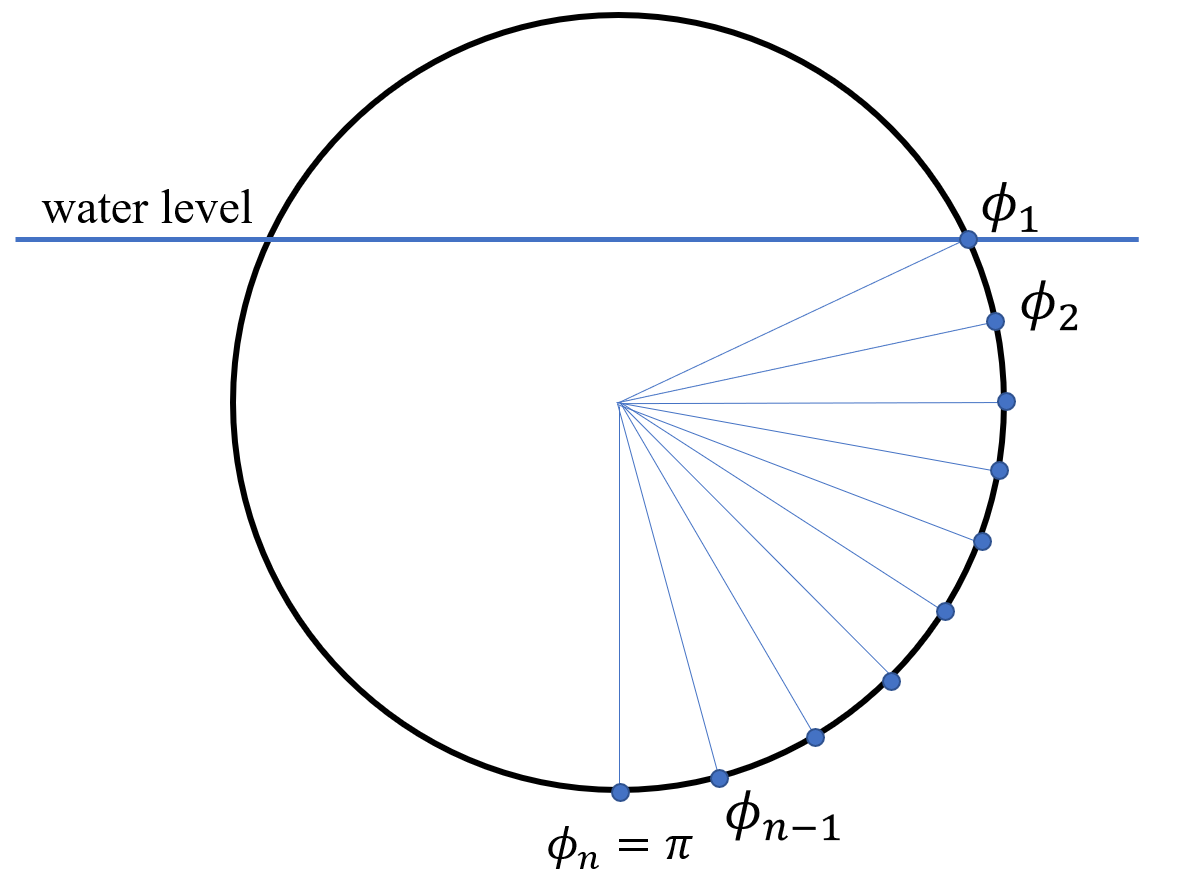}}} %
    \label{fig:reimanna}
    \qquad
    \subfloat[\centering Lower and Upper Sums]{{\includegraphics[scale=0.22]{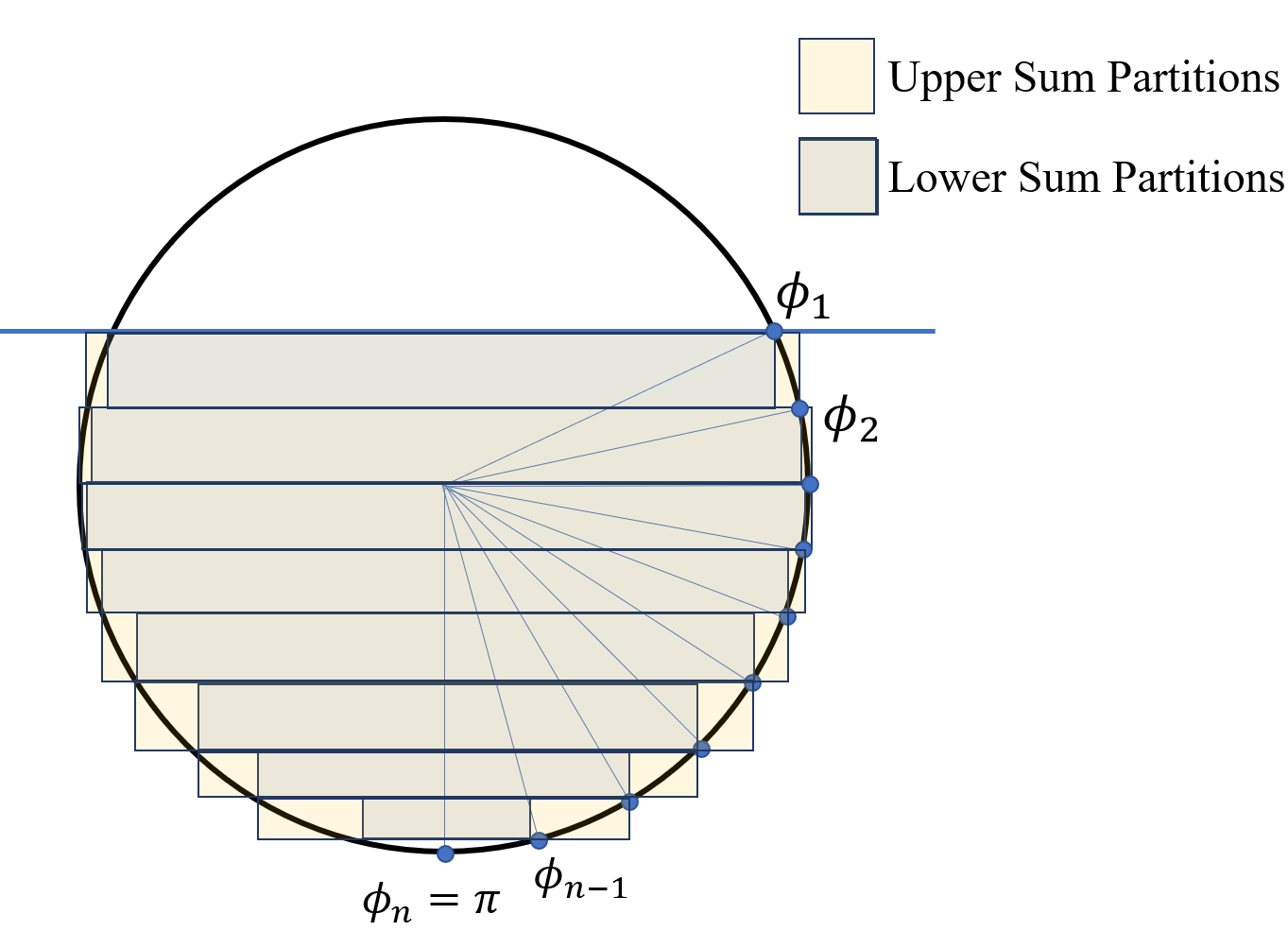}}} 
    \label{fig:reimannb}
    \caption{Discretization of the Submerged Volume}%
    \label{fig:reimann}%
\end{figure}

\begin{figure}
\centering
\includegraphics[scale=0.35]{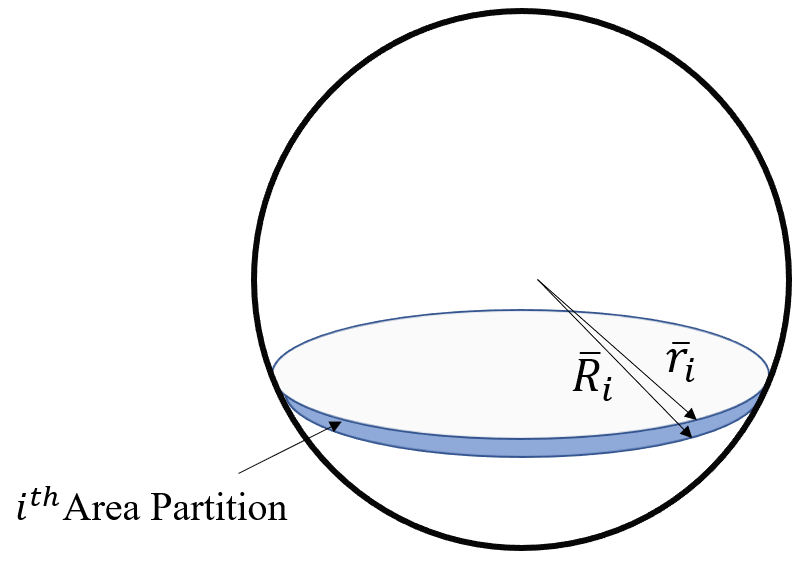}
\caption{$i^{th}$ Area Partition}
\label{fig:reimann2}
\end{figure}

The submerged volumes of the VSB and the FSB are calculated using the Riemann integrals approach, \cite{heywood_1965} where the submerged volume is divided into a set of $n$ horizontal disks (partitions), and the total submerged volume of the buoy is the sum of volumes of disks below the water surface, as shown in Fig.~\ref{fig:reimann}. 
Let $\phi_1$ be the angle of the highest wet disk on the VSB surface as shown in Fig.~\ref{fig:reimann}-\textcolor{blue}{a}; this can be written as: \begin{equation}
    \phi_1 = \phi(\max \{r_{dms_3}^{wet}\})
\end{equation} Recall that $r_{dms_3}$ is the vertical component (in the $s_3$ direction) of the $\vec{\boldsymbol{r}}_{dms}$ vector. Then it is clear that there exists a closed and bounded set that divides the circumference of the buoy into a set of partitions such that:
\begin{align}
   r_{dms_3}(\phi_n = \pi)   \leqslant  r_{dms_3}(&\phi_{n-1})    \leqslant \ldots \nonumber \\ & \leqslant r_{dms_3}(\phi_2) \leqslant r_{dms_3}(\phi_1)
\end{align} The height of the $i^{th}$ disk is calculated as $\Delta {r_{dms_3}}_{i} = {r_{dms_3}}_{i} -  {r_{dms_3}}_{i-1}$. Let $\bar{R}_i$ be the supremum (sup) of $r_{dms_1}$ in the $i^{th}$ disk, that is:
\begin{equation}
    {\bar{R}_i} = \sup_{\phi \in[\phi_{i-1},\phi_{i}]} \left(r_{dms_1}(\phi)\right).
\end{equation}
Likewise, let $\bar{r}_i$ be the infimum (inf) of $r_{dms_3}$ in the $i^{th}$ disk, that is:
\begin{equation}
   {\bar{r}_i} = \inf_{\phi \in[\phi_{i-1},\phi_{i}]} \left(r_{dms_1}(\phi) \right).
\end{equation}


Consider the volume of the submerged disks ($V_s$), there is a lower and upper limit for the volume as shown in Fig.~\ref{fig:reimann}-b. 
It is possible to compute the lower Riemann sum, $L(r_{dms_3},V_s)$, and upper Riemann sum, $U(r_{dms_3},V_s)$, for the submerged volume as follows:
\begin{align} 
U(r_{dms_3},V_s) &= \sum^{n}_{i=1} \pi {\bar{R}_i^2} \Delta {r_{dms_3}}_{i} \label{eq:73}\\
L(r_{dms_3},V_s) &= \sum^{n}_{i=1} \pi {\bar{r}_i}^2 \Delta {r_{dms_3}}_{i} \label{eq:74}
\end{align}
The difference between $U(r_{dms_3},V_s)$ and $L(r_{dms_3},V_s)$ is bounded; that is:
\begin{equation}
    U(r_{dms_3},V_s) - L(r_{dms_3},V_s) < \varepsilon_{V_s}
\end{equation}
where $\varepsilon_{V_s} > 0 $. An accurate calculation of the submerged volume would have a small $\varepsilon_{V_s}$. Clearly, as $n\rightarrow \infty$, the $\varepsilon_{V_s} \rightarrow 0$. However, as $n$ increases the computational cost increases significantly.

\subsubsection{Submerged Area Calculation} \label{subsec:sub_area_calc}

The internal surface area of the buoy was calculated using Riemann sums as well. The internal area is divided into $n$ number of horizontal slices with infinitesimal heights. Fig. (\ref{fig:reimann2}) shows a schematic for the $i^{th}$ area partition. The supremum $\bar{R}_i$and infimum $\bar{r}_i$ of the radii of this $i^{th}$ infinitesimal partition are calculated using Eqs. \eqref{eq:86} and \eqref{eq:87}  

\begin{align} 
    \bar{R}_i &= \sup_{\phi \in[\phi_{i-1},\phi_{i}]} \norm{r_{dms}(\phi)} \label{eq:86}\\
    \bar{r}_i &= \inf_{\phi \in[\phi_{i-1},\phi_{i}]} \norm{r_{dms}(\phi)} \label{eq:87}
\end{align} The height of the $i^{th}$ partition ($\Delta {r_{dms_3}}_{i} = {r_{dms_3}}_{i} -  {r_{dms_3}}_{i-1}$) and the lower and upper surface Riemann sums are calculated as follows: \begin{align} 
U(f,A_i) &= \sum^{n}_{i=1} 2 \pi \bar{R}^2_{i} \Delta {r_{dms_3}}_{i} \label{eq:88}\\
L(f,A_i) &= \sum^{n}_{i=1} 2 \pi \bar{r}^2_{i} \Delta {r_{dms_3}}_{i} \label{eq:89} 
\end{align} $\exists$ $\varepsilon_A > 0 $ such that
\begin{equation}
    U(f,A_i) - L(f,A_i) < \varepsilon_A
\end{equation} To increase the accuracy of the discretized area calculation, the number of partitions $n$ has to be big enough such that  $\varepsilon_{A} \rightarrow 0$ as $n\rightarrow \infty$. The final surface area of any portion is calculated as:

\begin{equation} \label{eq:114_5}
A_{i} =  \frac{ U(r_{dms_3},{A}_i) + L(r_{dms_3},{A}_i)}{2} 
\end{equation}

Noting that, using either of the areas calculated in Eq. \eqref{eq:88} and \eqref{eq:89} produce a first order accurate area calculation, on the other hand, Eq. \eqref{eq:114_5} produces a second-order accurate area calculation.

\subsubsection{Generalized Hydro Force}
The hydrodynamic force at the $i^{th}$ disk on the buoy surface can be expressed as: 
\begin{align}
    \vec{f}_i^{hydro} &= - {F_{hydro}}_i   {\hat{\boldsymbol{c}}_3}_i \label{eq:57} \\
   & =- \left(  {f_{hydrod}}_i \cos{\omega t}   +  \frac{ {(\rho_w V_s g)}_i }{ \cos (\pi - \psi_i)}\right) {\hat{\boldsymbol{c}}_3}_i
\end{align}
where ${(\rho_w V_s g)}_i$ is the buoyant force on the $i^{th}$ disk of the submerged volume.
The generalized hydro force can be written in the following form:
\begin{equation} \label{eq:60_1z} 
  \mathbit{Q}_j^{hydro} = \vec{\boldsymbol{f}}^{hydro} \cdot \frac{\partial \vec{r}_{dma}}{\partial q_j} , \, j = 1,2, ..., 6+N
\end{equation}

We can then write the hydro generalized forces using Eq. \eqref{eq:57} and Eq. \eqref{eq:60_1z} as follows:

\begin{align*}
\mathbit{Q}_{1,i}^{hydro} &=  \vec{\boldsymbol{f}}^{hydro}_i   \cdot \frac{\partial  \vec{\boldsymbol{r}}_{dma} }{\partial \vec{\boldsymbol{r}}_{sa,1}} =  - {F_{hydro}}_i  {\hat{\boldsymbol{c}}_3}_i \cdot \boldsymbol{1}_1 \nonumber \\ &= - {F_{hydro}}_i  \boldsymbol{1}_1^T {\hat{\boldsymbol{c}}_3}_i  = - {F_{hydro}}_i \sin{\psi_i}  \\
 \mathbit{Q}_{2,i}^{hydro}&=  \vec{\boldsymbol{f}}^{hydro}_i   \cdot \frac{\partial  \vec{\boldsymbol{r}}_{dma} }{\partial \vec{\boldsymbol{r}}_{sa,2}} =- {F_{hydro}}_i {\hat{\boldsymbol{c}}_3}_i \cdot \boldsymbol{1}_2 =0 \\
 \mathbit{Q}_{3,i}^{hydro} &= \vec{\boldsymbol{f}}^{hydro}_i   \cdot \frac{\partial  \vec{\boldsymbol{r}}_{dma} }{\partial \vec{\boldsymbol{r}}_{sa,3}}  = - {F_{hydro}}_i {\hat{\boldsymbol{c}}_3}_i \cdot \boldsymbol{1}_3  \nonumber \\ &=  - {F_{hydro}}_i \boldsymbol{1}_3^T {\hat{\boldsymbol{c}}_3}_i =- {F_{hydro}}_i \cos{\psi_i}    \\
 \mathbit{Q}_{4,i}^{hydro} &= \vec{\boldsymbol{f}}^{hydro}_i \cdot \frac{\partial  \vec{\boldsymbol{r}}_{dma} }{\partial \theta_1 } = 0 , \:\:\:
\hdots, \nonumber \\
 \mathbit{Q}_{6,i}^{hydro} &= \vec{\boldsymbol{f}}^{hydro}_i \cdot \frac{\partial  \vec{\boldsymbol{r}}_{dma} }{\partial {\theta_3} } = 0 \nonumber \\
 \mathbit{Q}_{7,i}^{hydro} &= \vec{\boldsymbol{f}}^{hydro}_i \cdot \frac{\partial  \vec{\boldsymbol{r}}_{dma} }{\partial {\eta_1} } = - {F_{hydro}}_i {\hat{\boldsymbol{c}}_3}_i  \cdot \boldsymbol{\Phi}_e(:,1) , \:\:\:
\hdots, \nonumber \\
 \mathbit{Q}_{6+N,i}^{hydro} &=  \vec{\boldsymbol{f}}^{hydro}_i \cdot \frac{\partial  \vec{\boldsymbol{r}}_{sa} }{\partial {\eta_N} } = - {F_{hydro}}_i {\hat{\boldsymbol{c}}_3}_i  \cdot \boldsymbol{\Phi}_e(:,N) 
\end{align*} 

The generalized hydro force on the buoy's shell is then expressed as:

\begin{align}
 &   \mathbit{Q}^{hydro} = \nonumber \\  &- \sum_{i=1}^{m} \begin{bmatrix}
    [0  & 0 &  - {F_{hydro}}_i \cos{\psi_i}]^T   
    \\ &\boldsymbol{0}& \\ 
    [F_{{hydro}_i} {\hat{\boldsymbol{c}}_3}_i  \cdot \boldsymbol{\Phi}_e(:,1) & 
     \hdots  
    &F_{{hydro}_i} {\hat{\boldsymbol{c}}_3}_i  \cdot \boldsymbol{\Phi}_e(:,N)]^T  \\
    \end{bmatrix} \label{eq:81}
\end{align}
where $m$ is the number of partitions on the buoy's shell.

\section{Model Validation}
The validation of this dynamic model follows the same approach as in \cite{robert2015modeling} .The natural frequencies for the VSB WEC surface obtained numerically via the Rayleigh–Ritz method are compared to the theoretical values obtained by solving Eqs.~\eqref{eq:7} and  \eqref{eq:11}, in Table~\ref{tab:1}. It is observed that the numbers are identical for the breathing mode ($n=0$), while for the other modes, there are negligible discrepancies.
Noting that the assumption used to obtain Eq. \eqref{eq:47_4} has no significant effect on the accuracy of the solution. 

\begin{table}
\centering
\caption{ Natural Frequencies Resulted From Rayleigh-Ritz Method and The Analytical Full Formulation }\label{tab:1}
\begin{tabular}{|c|c|c|}
\hline
\textbf{n} & \textbf{Rayleigh-Ritz }  & \textbf{Full Formulation } \\
\textbf{} & $rad/sec$ &  $rad/sec$ \\ \hline
\textbf{0} & 51.434449987 &  51.434449987 \\ \hline
\textbf{1} & 62.153777189 &  62.994144454   \\ \hline
\textbf{2} & 86.525263883 &  86.823947954   \\ \hline
\textbf{3} &115.855327968 & 115.941895367  \\ \hline
\textbf{4} &146.580164630 &  146.610087987   \\ \hline
\textbf{5} &177.820039783 & 177.832640732   \\ \hline
\textbf{6} &209.298426663 & 209.304856685   \\ \hline
\end{tabular}
\end{table}

\begin{figure*}[p]%
    \centering
    \subfloat[\centering $n=0$ and $\eta_1$ ]{{\includegraphics[scale=0.6]{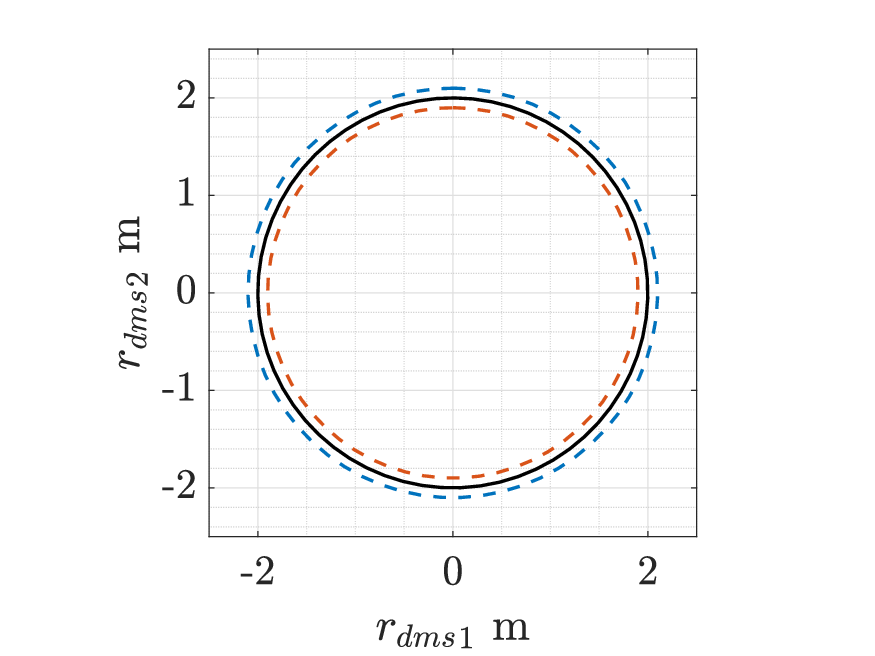}}} %
    \subfloat[\centering $n=1$ and $\eta_2$]{{\includegraphics[scale=0.6]{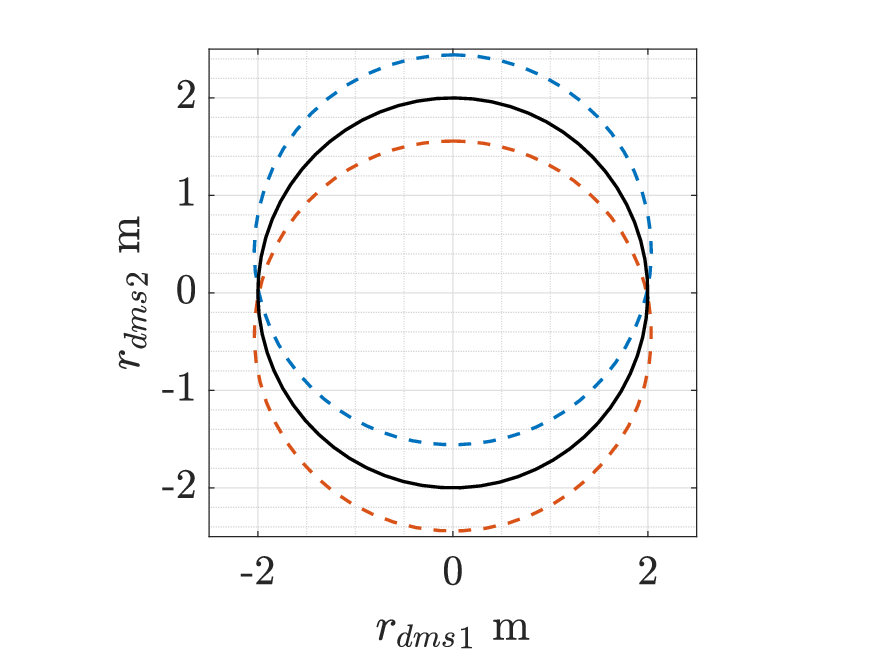} }} 
        \qquad
        \subfloat[\centering $\eta_3$ ]{{\includegraphics[scale=0.6]{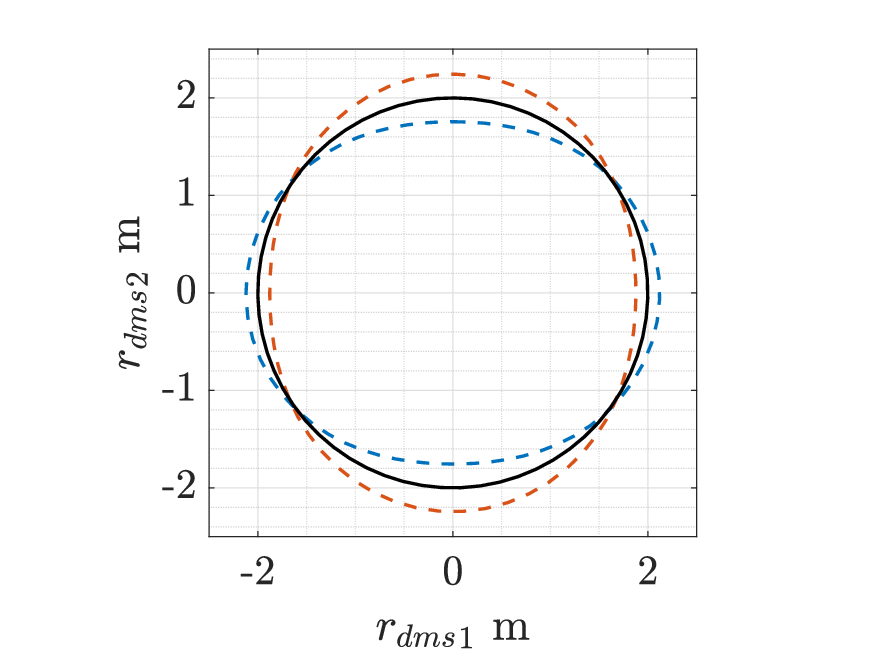}}} %
    \subfloat[\centering $\eta_4$]{{\includegraphics[scale=0.6]{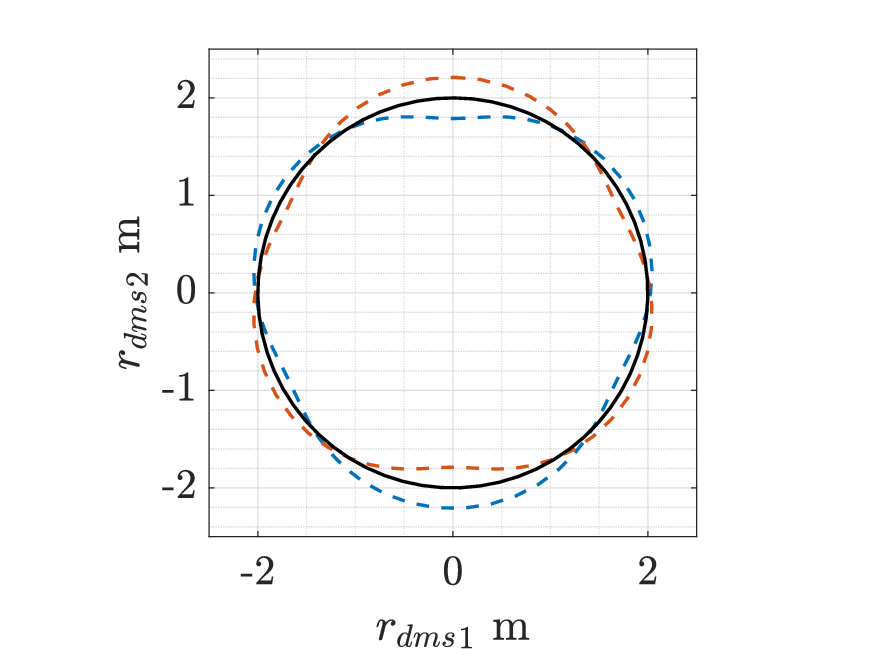} }}
            \qquad
            \subfloat[\centering $\eta_5$ ]{{\includegraphics[scale=0.6]{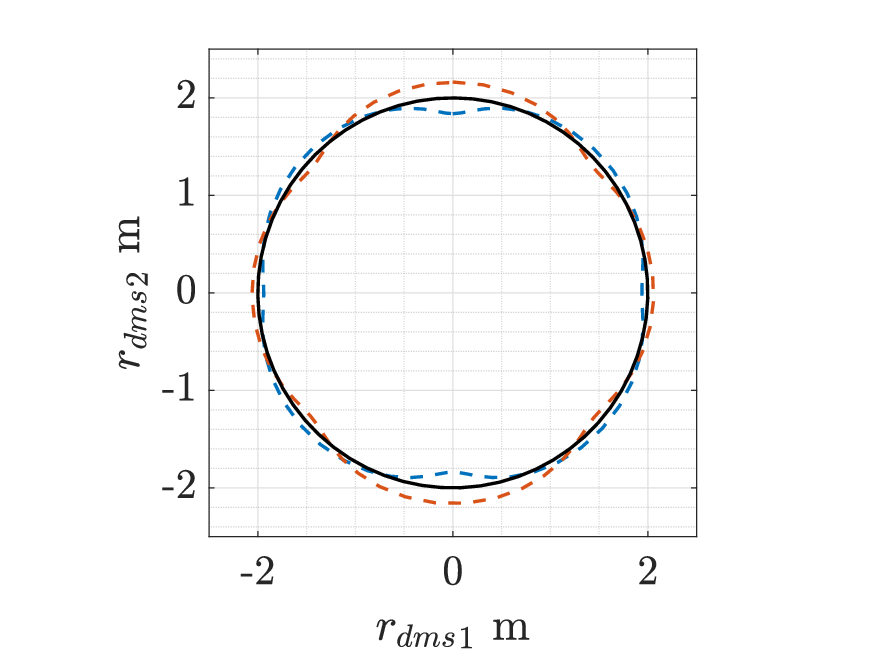}}} %
    \subfloat[\centering $\eta_6$]{{\includegraphics[scale=0.6]{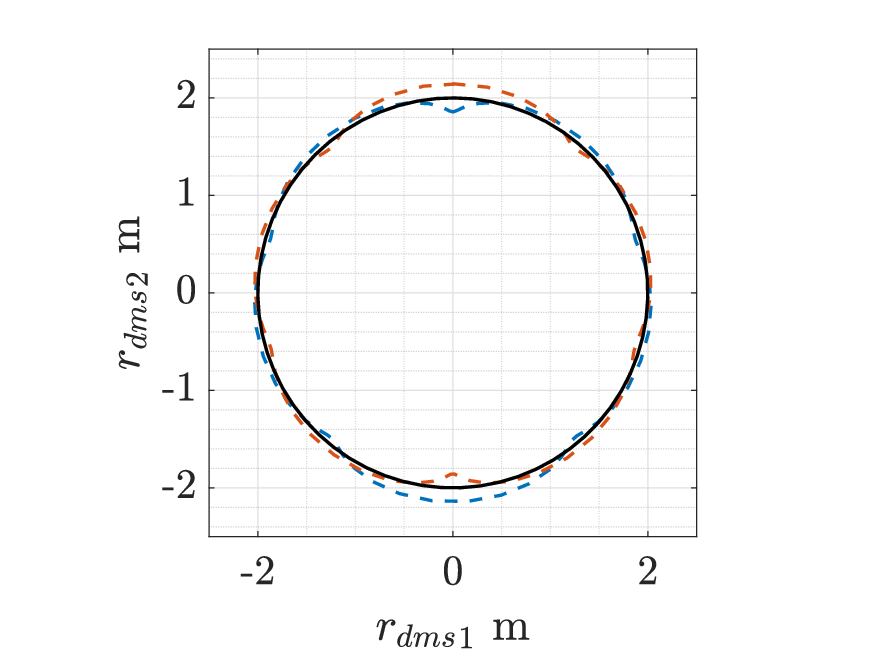} }} 
    \caption{Vibration Modes of axisymmetric spherical WECs, the black line is the undeformed buoy}%
    \label{fig:modes}%
\end{figure*}

Figure~\ref{fig:modes} shows the mode shapes resulted by the developed model. The breathing mode associated with $n=0$ and $\eta_1$ is presented in Fig.~\ref{fig:modes}-a where the buoy vibrate uniformly in the radial directions, Fig.~\ref{fig:modes}-b shows the mode shape associated with the rigid body motion ($\eta_1$), in which the buoy shell moves vertically. The third vibration mode shown in Fig.~\ref{fig:modes}-c transforms the shape of the buoy from prolate to oblate spheroids and vice versa. 

\section{Results and Discussion} \label{sec:results} 
Simulation results for the dynamic model of a spherical axisymmetric VSB WEC are presented in this section. The VSB is assumed to be made of a flexible hollow shell vented to the atmosphere. The simulations were carried for both the FSB WEC and VSB WEC; a free vibrating shell was tested as well as shells constrained from deforming at the top and middle horizontal section. The simulations also show the comparison between the energy converted by the VSB WEC and the FSB WEC. The equations of motion are solved using the MATLAB function $ode45$ which uses a six-step, fifth-order, Runge-Kutta method with variable time step.The simulation time for the tested VSB WECs is 2.4 hours while for the FSB WEC the simulation time is 2.11 minutes.

The radius of the buoy is $2$ m and the shell thickness is $0.01$ m. The modulus of elasticity and the poisons ratio are $10$ MPa and $0.3$, respectively. 
The total mass of the buoy is 17170 kg, and the wave excitation pressure is 1800 Pa with a period of $2.5$ sec. The damping coefficient of the PTO unit is set to $c = 8000$  Ns/m. All the initial conditions for the simulation were set to zeros except for $\vec{\dot{\boldsymbol{r}}}_{sa} =  [0 \:\: 0 \:\: -0.8]^T$ and the Legendre polynomial is truncated in 7 terms, i.e. $N=7$.


\begin{figure}[h]
\centering
\includegraphics[scale=0.6]{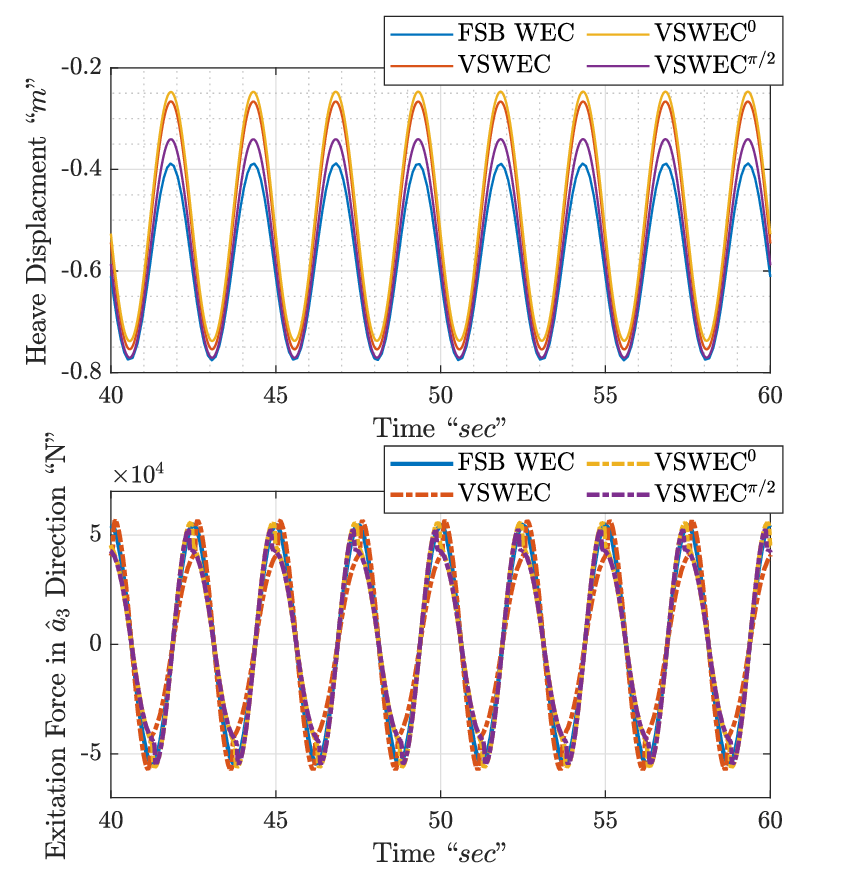}
\caption{Heave Displacement and The vertical Component of the Excitation Force}
\label{fig:res1}
\end{figure}

Figure~\ref{fig:res1} shows the heave displacement and the vertical component of the excitation forces for all the four cases. The transient effects die out after almost 25 seconds in the simulation, and the plots show the interval from 40 to 60 seconds. It is observed that the three VSB WECs designs have higher displacements compared to the FSB WEC. Larger motions are usually associated with higher energy conversion, and hence higher displacements are usually desirable. These results support the hypothesis of this research which is that flexible buoys would leverage the waves and behave like a rigid buoy that has reactive power; indeed, the reactive power in this case is obtained from the waves themselves.

\begin{figure}
    \centering
    \subfloat[\centering The Heave Velocity ]{{\includegraphics[scale=0.55]{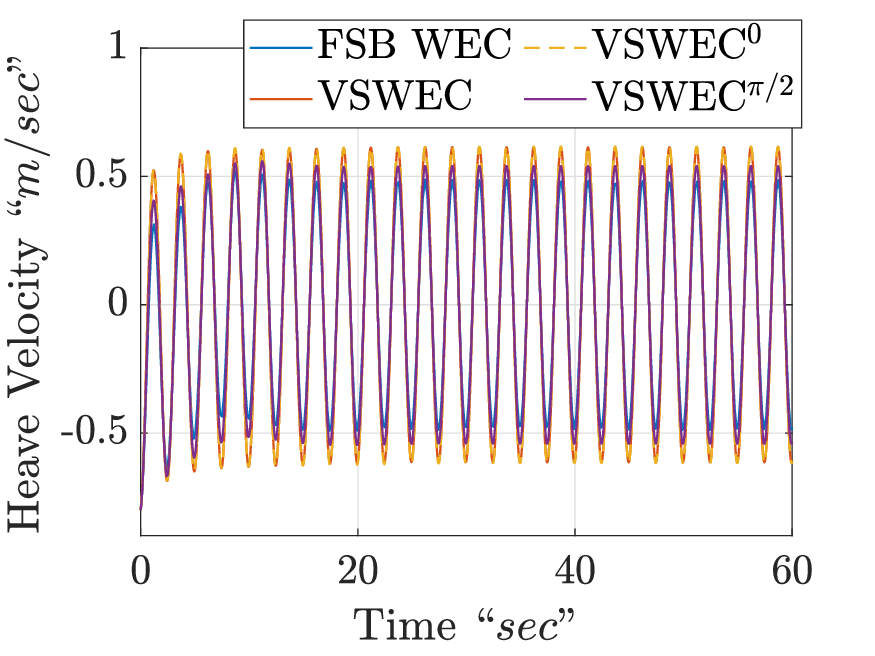}}} %
    \qquad
    \subfloat[\centering The PTO Force]{{\includegraphics[scale=0.55]{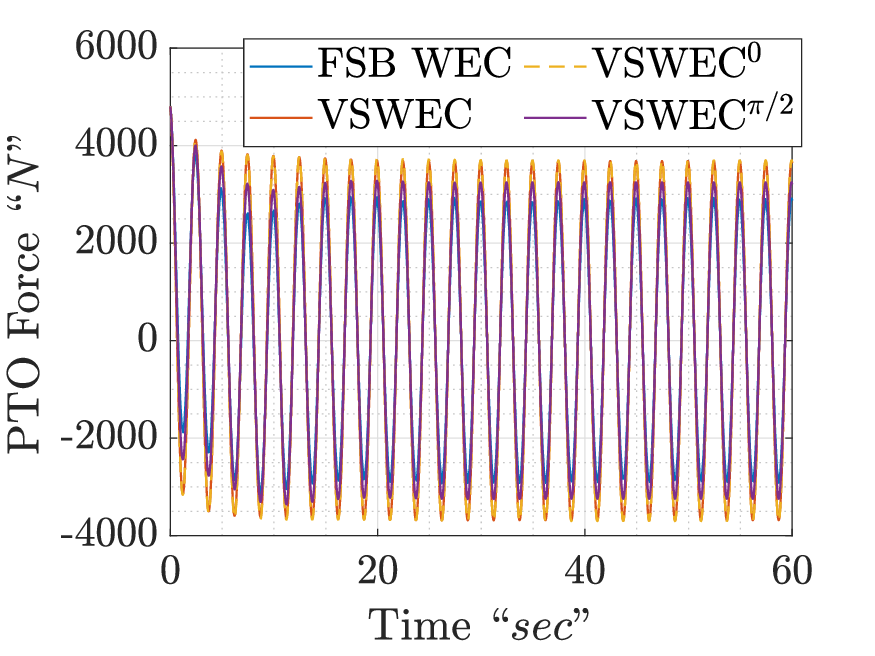} }} 
    \caption{Heaving Velocity and PTO Force}%
    \label{fig:res2}%
\end{figure}

\begin{figure}
\centering
\includegraphics[scale=0.55]{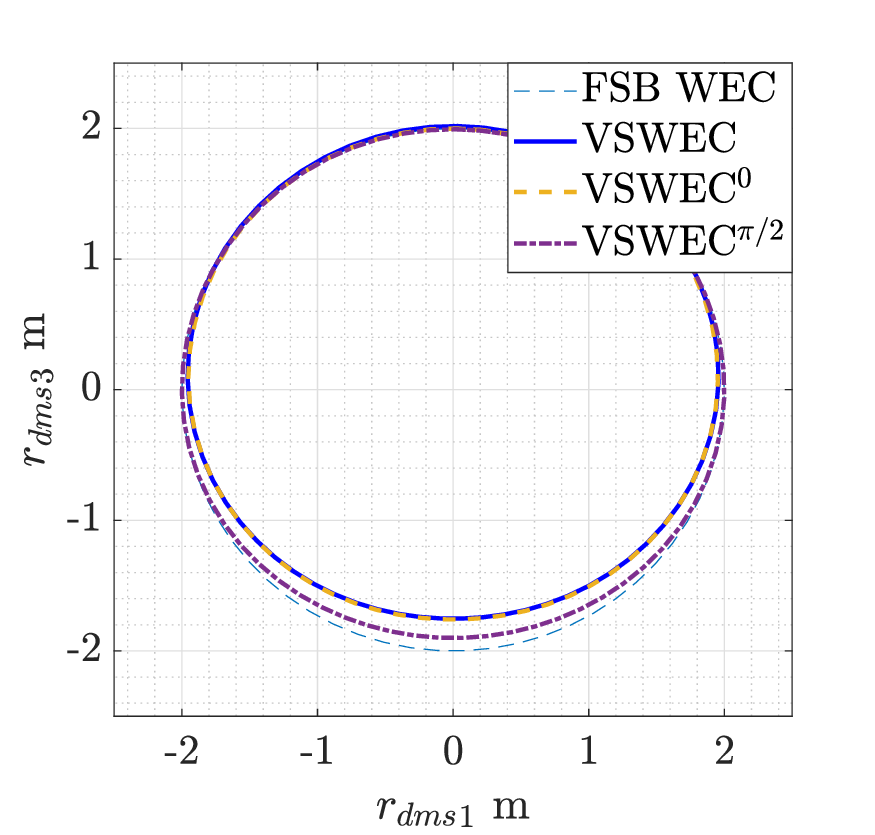}
\caption{Shapes of the FSB and the VSBs at their state of minimum vertical deformation}
\label{fig:res1_1}
\end{figure}

It is also noticed that the peak-to-peak (pk-pk) displacement of the VSWEC$^{0}$ is slightly higher than that of the VSWEC; this is due to the bigger deformations of the bottom half of the VSWEC$^{0}$ buoy. On the other hand, the peak-to-peak displacement of the VSWEC$^{\pi/2}$ is less than the other VSB WEC designs and higher than the FSB WEC because the former allows for a pivoting point at angle $\phi_c = 90^o$. 
Also, a phase shift of $\pi/4$ is noticed between the vertical component of the excitation force and the heave displacement occurred due to the high non-linearity of the VSB WEC.

Figure~\ref{fig:res1_1} shows the shapes of the FSB and the VSBs at their state of minimum vertical deformation (for highlighting the difference in deformations, the deformation vector is multiplied by a factor of 10), in the $-\hat{c}_3$ direction. Fig~\ref{fig:res1_1} demonstrates the imposed constraint on the VSWEC$^0$ since it coincides with the FSB (i.e. no deformation) at the top point at $\phi_c = 0^o$, which is the constraint point. Likewise, Fig.~\ref{fig:res1_1} shows that the VSWEC$^{\pi/2}$ coincides with the FSB WEC (i.e. no deformation) at $\phi_c = 90^o$  which is the location of the imposed no-deformation constraint. 

The heave velocities of the VSW WECs are higher than the FSB WEC as shown in Fig.~ \ref{fig:res2}-a. At steady state response, the waveform for the heave velocity for the VSWEC and the VSWEC$^{0}$ coincide over each other 
The PTO force shown in Fig. \ref{fig:res2}-b is calculated by multiplying $-c$ by the heave velocity as expressed in Eq.\eqref{eq:73_3}.

\begin{figure*}
\centering
\includegraphics[scale=0.7]{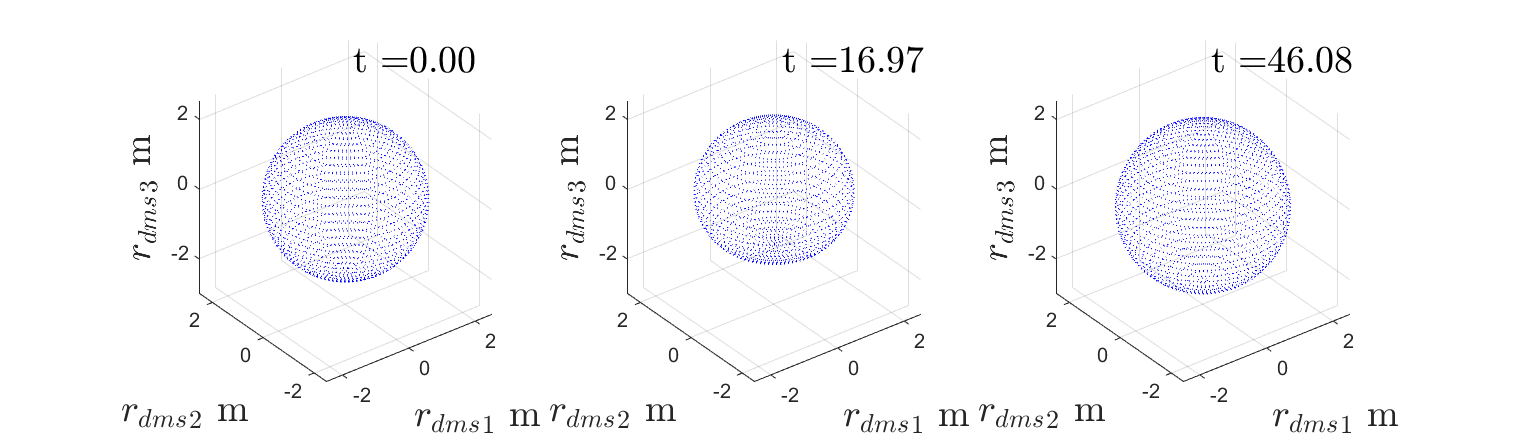}
\caption{3-D Visualization for the VSWEC in various time steps}
\label{fig:res1_11}
\end{figure*}

Figure~\ref{fig:res1_11} shows a 3D visualization VSB WEC at different time steps; a similar scaling factor was applied to the plot in Fig.~\ref{fig:res1_1}. At $t = 0$ the VSB is not deformed, i.e., the buoy shape is spherical.  At $t = 16.97$ and $46.08$ sec the VSBs have oblate and prolate spheroid shapes, respectively. At t = 16.97 sec the wave crest pushes the buoy out of the water, resulting in an oblate spheroid shape. This results in larger $\mathbit{Q}_{3,i}^{hydro}$ because of the \textit{cosine} term in Eq.~\ref{eq:81}, which is the component responsible for the vertical component of the hydro force acting on the CG. At t = 46.08 sec the buoy is encountering a wave trough and a prolate spheroid shape results in an opposite effect compared to the oblate spheroid and $\mathbit{Q}_{3,i}^{hydro}$ is reduced, i.e, the buoy dives more into the water. The switching between oblate and prolate spheroid shapes with the wave crests and troughs results in higher pk-pk displacements and velocities for the VSB WECs compared to the FSB WEC as shown in Figs.~\ref{fig:res1} and \ref{fig:res2}.

\begin{figure}
    \centering
    \subfloat[\centering WECs Total Volume ]{{\includegraphics[scale=0.58]{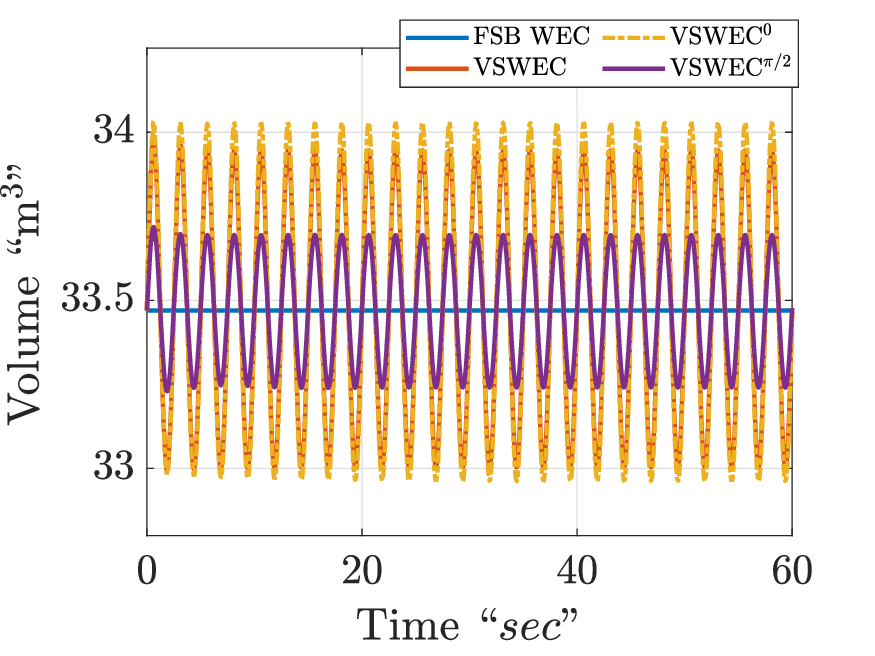}}} %
   \qquad
    \subfloat[\centering WECs Total Surface Area]{{\includegraphics[scale=0.58]{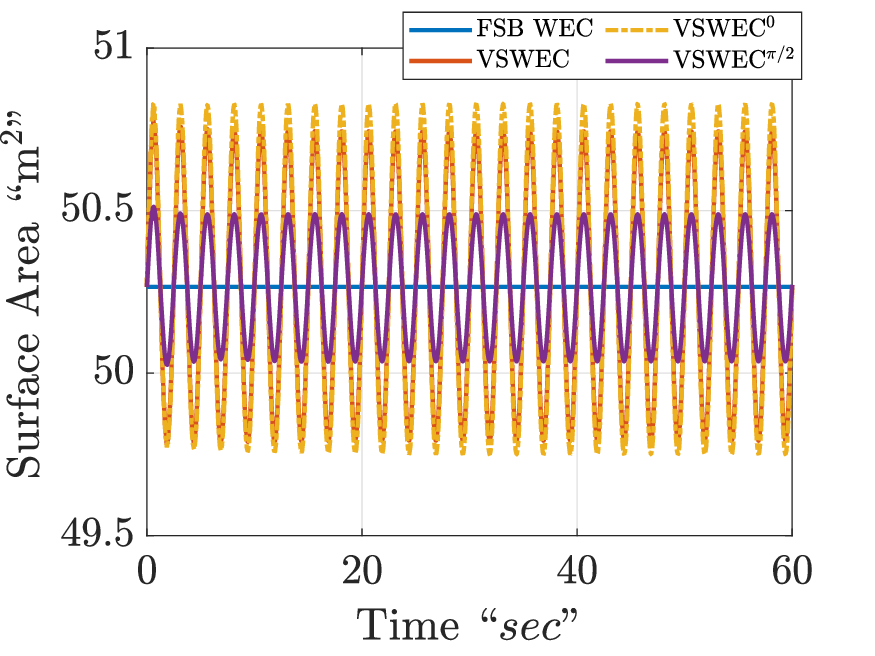} }} 
    \caption{Volume and surface area cover time}%
    \label{fig:area_and_vol}%
\end{figure}

The total volume and total surface area change over time are shown in Figures~\ref{fig:area_and_vol}-a and \ref{fig:area_and_vol}-b. The change in the volume and surface areas corresponds to the change in the buoyant force and excitation forces as discussed in subsection \ref{subsec:3.1}. It is noticed that the VSWEC$^{0}$ has the highest pk-pk change for the volume and areas change, and the VBWEC$^{\pi}$ has the least change in volume and surface area.

\begin{figure}
    \centering
    \subfloat[\centering Generated Power ]{{\includegraphics[scale=0.5]{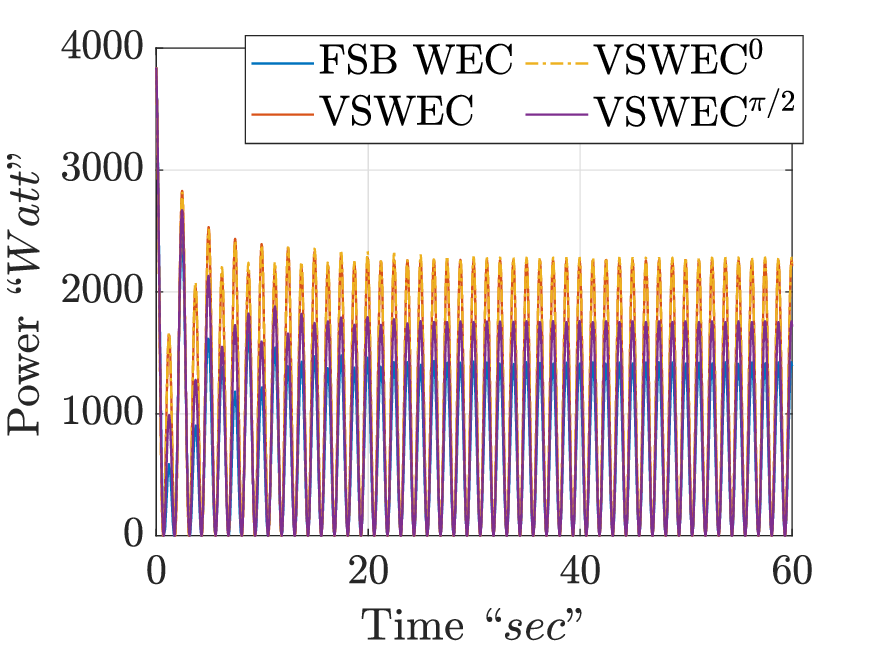}}} %
    \qquad
    \subfloat[\centering Harvested Energy]{{\includegraphics[scale=0.5]{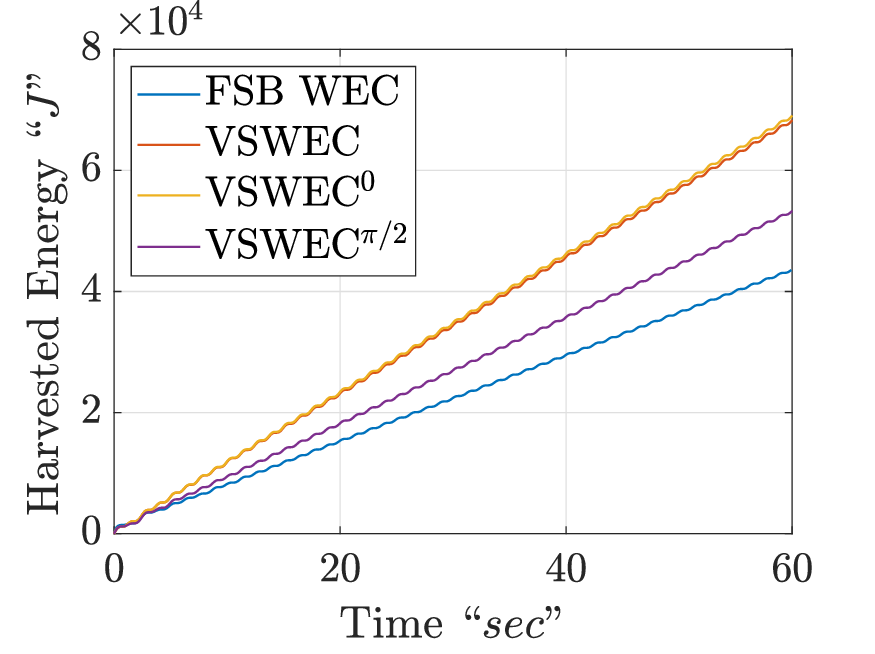}}} 
    \caption{Generated Power and Harvested Energy}%
    \label{fig:res3}%
\end{figure}

The generated power peaks for all the VSB WEC designs are higher than the FSB WEC as shown in Fig. \ref{fig:res3}-a; the lowest power peaks are generated by the FSB WEC. Moreover, the peaks of the VSWEC$^0$ and VSWEC are almost overlapping. Fig.~\ref{fig:res3}-b shows the total harvested energy over a period of 60 seconds for all four cases. Clearly, there is a multiple-fold increase in the harvested energy of a VSB WEC compared to a FSB WEC. The VSWEC$^{0}$ harvests 59$\%$ more energy with reference to the FSB WEC; the VSWEC and the VBWEC$^{\pi/2}$ harvested 57$\%$ and 22.8$\%$ more energy, respectively. The discussion of the effect of material properties on the VSB WEC is beyond the scope of the current article; however, it is worth noting that several values of the modulus of elasticity were tested, and it is observed that the softer the buoy material, the more energy is harvested from the waves. Finally, Table~\ref{tab:2} summarizes all the performance measures for the four test cases, where the generated energy is the energy harvested over 60 seconds of the simulation period.

\begin{table*}
\centering
\caption{Comparison Between the Overall Performance of the Different Designs for WECs at Steady State based on the assumed hydrodynamic pressures}
\begin{tabular}{c|c|c|c|c|c|}
\cline{2-6}
\textbf{} & \textbf{displacement} & \textbf{Velocity}& \multicolumn{1}{l|}{\textbf{Control Force}}& \multicolumn{1}{l|}{\textbf{Power}} & \multicolumn{1}{l|}{\textbf{Generated Energy}} \\
\multicolumn{1}{l|}{}& (m pk-pk)& \multicolumn{1}{l|}{(m/sec pk-pk)} & (N pk)& (Watt)& (KJ)\\ \hline
\multicolumn{1}{|c|}{\textbf{FSWEC}}& 1.1632& 0.9732 & 109480 & 1421 & 43.42\\ \hline
\multicolumn{1}{|c|}{\textbf{VSWEC$^{\pi/2}$}}& 1.1122 & 1.081 & 104360 & 1757 & 53.31\\ \hline
\multicolumn{1}{|c|}{\textbf{VSWEC}}& 1.019 & 1.2266 & 113160 & 2253 & 68.27\\ \hline
\multicolumn{1}{|c|}{\textbf{VSWEC$^0$}} & 0.9844 & 1.2338 & 111220 & 2283 & 69.05\\ \hline
\end{tabular}\label{tab:2}
\end{table*}


\section{Conclusion and Future Work} \label{sec:conclusion}
The equations of motion of spherical Variable-Shape Buoy Wave Energy Converters, with axisymmetric deformations, were derived using a Lagrangian formulation in this paper. A Rayleigh-Ritz method along with the classical bending theory for stress-strain relations were used to approximate the equations of motion to finite-dimension equations of motion, in the six degrees of freedom. Holonomic constraints were imposed to limit the buoy to only-heave motion, and to enforce no-deformation at specific locations, to account for the power take-off unit installation flanges. The inner volume of the VSB WEC is assumed to be vented to the atmosphere to exclude any internal pressure variation effect on the buoys' shell. 
The numerical results support the hypothesis of this work which is that a VSB WEC would harvest energy at a significantly higher rate compared to that of a FSB WEC, when both WECs use no reactive power. 
The VSB WEC with zero deformation at its highest vertical location harvested more energy than the VSB WEC with unconstrained shell deformations and the VSB WEC with shell constraint at the horizontal midsection.

For future work, one-way and two-way FSI investigations are required to examine the hydrodynamic performance of the VSB WECs in realistic regular and irregular waves environment; in these FSI model, the developed dynamic model would be coupled with a fluid solver (e.g., OpenFoam, Nemoh, WAMIT) such that at every time step the deformed geometry is exported to the fluid solver where the pressure distribution around the VSB WEC shell is calculated. In case of using BEM solvers, expressions for the generalized added mass, damping, hydostatic stiffness and excitation force coefficients can be calculated in the generalized coordinates. 


\section{Acknowledgment}
This material is based upon work supported by the National Science Foundation (NSF), USA, under Grant Number 2023436.
The research reported in this paper is partially supported by the HPC@ISU equipment at Iowa State University, some of which has been purchased through funding provided by NSF under MRI grant number 1726447.

%

\bibliographystyle{IEEEtran}
\bibliography{references}

\end{document}